\documentclass[a4paper,oneside,12pt]{amsart}
\usepackage{etex}

\usepackage{AJheaderArticle}

\title[Integration on $1$-currents]{A Henstock-Kurzweil type integral on $1$-dimensional integral currents}
\author{Antoine Julia}
\date{\today}
\address{Dipartimento di Matematica, via Trieste 63, 35121 Padova, Italy}
\email{antoine.julia@imj-prg.fr} 

\begin{document}
\begin{abstract}
We define a non-absolutely convergent integration on integral currents of dimension $1$ in Euclidean space. This integral is closely related to the Henstock-Kurzweil and Pfeffer Integrals. Using it, we prove a generalized Fundamental Theorem of Calculus on these currents. A detailed presentation of Henstock-Kurzweil Integration is given in order to make the paper accessible to non-specialists.
\end{abstract}
\maketitle
\section{Introduction}
The goal of this paper is to give a comprehensive presentation of an integration method for functions defined on the support of an integral current of dimension $1$ in Euclidean spaces. This method is inspired from the Henstock-Kurzweil (HK) and Pfeffer Integrals \cite{Kurzweil1957generalized, Henstock1961definitions, pfeffer1991gauss}, and tailored for the study of the Fundamental Theorem of Calculus. The HK Integral is a variant of the Riemann Integral, yet it is more general than the Lebesgue Integral --- all Lebesgue integrable functions are HK integrable --- but {\it non absolutely convergent}: there exist functions which are HK integrable, while their absolute value is not; in the same way that the series $\sum_k (-1)^kk^{-1}$ converges, while $\sum_k k^{-1}$ does not. For functions defined on a bounded interval $[\,a,b\,]$, the Fundamental Theorem of Calculus of HK integration is
\begin{Theorem}\label{fti1}
  Let $f: [\,a,b\,] \to \R$ be a continuous function which is differentiable everywhere, then its derivative $f'$ is HK integrable on $[\,a,b\,]$ and there holds:
  \begin{equation}
    \label{eq:fti1}
    (HK)\int_a^b f' = f(b)-f(a).
  \end{equation}
\end{Theorem}
Note that some other integration methods have been defined which also satisfy Theorem \ref{fti1} in particular a ``minimal'' theory in \cite{BonDipPre2000constructive}. It is also interesting to note that a small variation in the definition of the HK integral yields the Mac Shane Integral \cite{McShane1973Unified} which is equivalent to the Lebesgue Integral on an interval.

The Riemann-like formulation of the HK integral makes it straightforward to allow for singularities in the above theorem: if $f$ is only differentiable at all but countably many points of $[\,a,b\,]$, the result still holds. This statement is in some sense optimal. Indeed, as shown by Z. Zahorsky in \cite{Zahorski1946ensemble}, the set of non-differentiability points of a continuous function is a countable union of $G_\delta$ sets. In particular, if it is uncountable, it must contain a Cantor subset by \cite[Lemma 5.1]{Oxtoby1980MnC}. Finally, to any Cantor subset of an  interval having zero Lebesgue measure, one can associate a ``Devil's Staircase'' which has derivative equal to $0$ almost everywhere and is non constant. 

However, the differentiability condition can be relaxed and replaced by a pointwise Lipschitz condition. Thus a more general statement is
\begin{Theorem}\label{fti2}
  Let  $f: [\,a,b\,] \to \R$ be a continuous function which is pointwise Lipschitz at all but countably many points, then, it is differentiable almost everywhere in $[\,a,b\,]$, its derivative $f'$ is HK integrable on $[\,a,b\,]$ and identity \eqref{eq:fti1} holds.
\end{Theorem}
 
Natural extensions of the Fundamental Theorem of Calculus include the Gauss Green (or Divergence) Theorem and Stokes' Theorem. For the former in bounded sets of finite perimeter an integral has been developped by W. F. Pfeffer in \cite{pfeffer1991gauss}, after works of J.  Mawhin \cite{Mawhin1981generalized} and J. Ma\v{r}\`ik \cite{Marik1965extensions}. The results extend naturally to Stokes' Theorem on smooth oriented manifolds.

  For singular varieties, an integral adapted to Stokes' Theorem has been defined by the author on certain types of integral currents in Euclidean spaces \cite{ThesisAJ,juliaARXIVstokes}. We present here the content of the second chapter of the author's thesis, where we focus on one dimensional integral currents. These are treated in a different way as they can be decomposed into a countable family of curves. We thus define an integral closer to the Henstock-Kurzweil one, which we call the HKP Integral. Given a current $T$, the set $\textrm{Indec}(T)$ consists of all the points of $\spt T$ which are in the support of an indecomposable piece of $T$ (see Sections \ref{section pieces} and \ref{sec 1d integration} for the notations). Our main result is the following:
\begin{restatable}[Fundamental Theorem of Calculus]{Theorem}{hkpfti}\label{FTI 1d}
Let $T$ be a fixed integral current of dimension $1$ in $\Rn$, and $u$ be a continuous function on $\spt T$. Suppose that $u$ is pointwise Lipschitz at all but countably many points in $\textrm{Indec}(T)$ and that $u$ is differentiable $\Vert T\Vert$ almost everywhere, then $x\mapsto \langle \dD u (x) , \vect T(x) \rangle$ is $\PfefferOne$ integrable on $T$ and
\[
     (\partial T)(u) = (\PfefferOne) \int_T \langle \dD u , \vect T \rangle.
\]
\end{restatable}
This theorem is equivalent to Theorem \ref{fti2} when $T$ represents an interval.

\subsection*{Summary of the paper}
In Section \ref{HK section}, we define the Integral of Henstock and Kurzweil and its main properties along with schemes of proofs of the main theorems. We also give an equivalent definition of integrability --- inspired from the Pfeffer Integral --- which will be useful in the sequel. It is important to note that the Pfeffer Integral is not equivalent to the HK Integral.

 In Section \ref{section pieces}, we recall the definition of integral currents of dimension $1$ in Euclidean spaces and define the main ingredients of HKP integration: pieces of a current and functions on the space of pieces of a current, we also study the derivation of these functions  (following Federer \cite[2.9]{FedererGMT}). Section \ref{sec 1d integration} contains the definition of HKP integration and the proof of its main properties, as well as the proof of Theorem \ref{FTI 1d}.

\subsection*{Possible generalizations}\label{section generalizations}
First we could ask if $u$ can be allowed to be discontinuous (yet bounded) outside of the density set of $\Vert T\Vert$. Proposition \ref{example circles} and Example \ref{example two curves} show that this is not straighforward.

A natural question would be whether Theorem \ref{FTI 1d} could be generalized to normal currents in Euclidean spaces. Indeed, normal currents of dimension $1$ also admit a decomposition into lipschitz curves. More precisely, by a Theorem of S. K. Smirnov \cite{Smirnov1993Decomp}, a current $T$ of dimension $1$, with finite mass and finite boundary mass in $\Rn$ can be written
\begin{equation}\label{eq: normal decomp}
    T= \int \lseg \gamma \rseg \dd \mu(\gamma),
  \end{equation}
  where $\mu$ is a finite measure on the space of Lipschitz curves. However, there is no a priori constraint on the measure $\mu$: it can certainly have a higher dimensional behavior. It is therefore impossible to work with countable sums of pieces and one would probably need another notion of piece of a normal current to follow the same plan as here. Recall that Fubini-type arguments do not work well with non-absolutely convergent integrals, as shown in \cite[Section 11.1]{Pfeffer1993riemannbook}, which indicates at the very least that one should be careful here. Note also that the space of curves, on which we would have to integrate is far from Euclidean.

Another natural idea would be to consider integral currents of dimension $1$ in Banach spaces or complete metric spaces, following \cite{AmbrKirch2000} or \cite{DPHardt2012rectifnflat}. The same strategy should work overall, although I do not know if the result can be attained in the same generality.

Finally I would like to mention that there are works on integration on more fractal objects with different methods \cite{Young1934limits, HarrisonNorton1991fractalcurves, Zust2011integration}.

\subsection*{Acknowledgements}
I wrote my PhD thesis at the Institut de Mathématiques de Jussieu, Université Paris Diderot USPC. I am currently supported by University of Padova STARS Project "Sub-Rieman\-nian Geometry and Geometric Measure Theory Issues: Old and New" (SUGGESTION). I wish to thank my PhD advisor Thierry De Pauw and my academic older brother Laurent Moonens, for their help during this thesis, as well Marianna Csörnyei for a helpful conversation on the Besicovitch Covering Theorem. Finally, wish to mention that the notions of pieces and subcurrents studied in this paper and in \cite{ThesisAJ, juliaARXIVstokes} are very close to the subcurrents defined by E. Paolini and E. Stepanov in \cite{PaoSte2012acyclic} for normal currents in metric spaces.

\section{The integral of Kurzweil and Henstock}\label{HK section}

\subsection{Definition and classical properties}

A nonnegative function defined on a set $E\subseteq \R$ is called a \emph{gauge} if its zero set is countable. In the classical definition of the Henstock-Kurzweil Integral, gauges are always positive, but for our purpose it makes sense to allow the gauge to take the value zero in a countable set.  A \emph{tagged family} in an interval $[\,a,b\,]$ is a finite collection  of pairs $([\,a_{j}, b_j\,],x_j)_{j=1,2,\dots,p}$ where one has $a\leq a_1<b_1\leq a_2<\dots\leq a_p <b_p \leq b$ and for all $j$, $x_j\in [\,a_{j-1},a_j \,]$ The \emph{body} of a family $\calP$ is the union denoted by $[\calP]$ of all the intervals in $\calP$. A \emph{tagged partition} in $[\,a,b\,]$ is a tagged family whose body is $[\,a,b\,]$.  If $\delta$ is a gauge on $[\,a,b\,] $, we say that a tagged family (or a tagged partition) is $\delta$-fine, when for all $j$, $b_j-a_{j}<\delta(x_j)$. In particular, there holds $\delta(x_j)>0$, for all $j$. 

\begin{Definition}\label{HK integral}
A function $f$ defined on a compact interval $[\,a,b\,]$ is \emph{Henstock-Kurzweil integrable on $[\, a,b\,]$} if there exists a real number $\alpha$ such that for all $\epsilon>0$, there exists a positive gauge $\delta$ on $[\,a,b\,]$ such that for each $\delta$-fine tagged partition $\calP =\{([\, a_{j-1},a_{j}\,],x_j)\}_{j=1,\dots,p}$, there holds:
\[
   \left\vert \sum_{j=1}^pf(x_j)(a_j-a_{j-1})-\alpha\right \vert <\epsilon.
\]  
\end{Definition}
In the following, we will write $\sigma(f,\calP)$ for the sum on the left hand side, whenever $\calP$ is a tagged family. If $\alpha$ as above exists, we denote it by $(HK)\int_a^bf$. This definition is well posed as a consequence of the following key result.
\begin{Lemma}[Cousin's Lemma]\label{HK Cousin}
    If $I$ is a closed bounded interval and $\delta$ is a positive gauge on $I$, then a $\delta$-fine tagged partition of $I$ exists.
\end{Lemma}
\begin{proof}
  Suppose no $\delta$-fine tagged partition of $I$ exists. Consider the two halves of $I$: $I_1$ and $I_2$. Either $I_1$ or $I_2$ does not admit a $\delta$-fine tagged partition. By successive divisions, we can find a decreasing sequence of closed intervals of the form $I^p = I_{j_1,j_2,\dots,j_p}$ where $j_k \in \{1,2\}$ and $I_{j_1,\dots,j_p,1}$ and $I_{j_1,\dots,j_p,2}$ are the two halves of $I_{j_1,\dots,j_p}$. We can choose the intervals $I^p$ for  $p =1,2,\dots$ so that none of them admits a $\delta$-fine tagged partition. There exists $x\in I\cap \bigcap_{p=1}^\infty I^p$. Since $\delta$ is positive on $I$, $\delta(x)>0$ and as $\diam(I^p) = 2^{-p} \diam I$, there exists $p$ such that $\diam I^p< \delta(x)$. This implies that $((I^p,x))$ is a $\delta$-fine tagged partition of $I^p$, a contradiction.
\end{proof}
The following propositions list the main properties of the HK integral. 
\begin{Proposition}[Cauchy Criterion for integrability]\label{Cauchy HK}
  A function $f$ is HK integrable on the interval $[\,a,b\,]$, if and only if for each $\epsilon>0$ there exists a positive gauge $\delta$ on $[\,a,b\,]$ such that whenever $\calP_1$ and $\calP_2$ are $\delta$-fine tagged partitions of $[\,a,b\,]$, there holds
\[
   \left \vert \sigma(f,\calP_1) -\sigma(f,\calP_2)\right \vert <\epsilon.
\]
\end{Proposition}
\begin{Proposition}\label{HK props}
 Let $f$ be a Henstock-Kurzweil integrable function on the interval $[\,a,b\,]$:
  \begin{enumerate}[label = (\arabic*)]
    \item  If $g$ is HK integrable on $[\,a,b\,]$ and $\lambda$ is a real number, then $f+\lambda g$ is HK integrable and 
\[
     (HK) \int_a^b (f+\lambda g) = \left ((HK) \int_a^b f\right ) + \lambda \left ( (HK) \int_a^b g \right ).
\] \label{HK linear}
    \item If a function $g$ is equal to $f$ almost everywhere on $[\,a,b\,]$, then $g$ is also HK integrable and has the same integral.
    \item If $g$ is Lebesgue integrable, it is also HK integrable and the two integrals coincide.
    \item The restriction of $f$ to a subinterval $[\,c,d\,]\subseteq [\,a,b\,]$ is HK integrable on $[\,c,d\,]$.
    \item (Saks-Henstock Lemma) \label{HK Saks Henstock} For $\epsilon>0$ and $\delta$ a positive gauge corresponding to $\epsilon$ in the definition of integrability of $f$, given any tagged family $(([\,a_j,b_j\,],x_j))_{j=1}^p$ in $[\,a,b\,]$ there holds
\[
    \sum_{j= 1} ^p \left \vert f(x_j) (b_j-a_j) -(HK)\int_{a_j}^{b_j} f\right \vert < 2\epsilon.
\]
    \item The function $F: [\,a,b\,]\to \R; x\mapsto (HK) \int_a^x f$ is continuous it is called the \emph{indefinite HK integral of $F$}. Also, if $f$ is nonnegative, $F$ is nondecreasing.\label{HK primitive}
    \item The function $F$ above is differentiable almost everywhere with derivative equal to $f$.\label{HK diff primitive}
    \item $f$ is Lebesgue measurable. \label{HK measurable}
    \item $f$ is Lebesgue integrable if and only if $f$ and $\vert f\vert$ are HK integrable.\label{HK abs to Lebesgue}
  \end{enumerate}
\end{Proposition}
The proofs of these results can be found in any treaty on Henstock-Kurzweil Integration (see Chapter 9 of \cite{Gordon1994integrals}, the recent book \cite{Moonens2017integration} --- in French, or the exercises in the appendix H to \cite{CohnMT2013}). In Section \ref{sec 1d integration} we prove results comparable to Proposition \ref{Cauchy HK} and Proposition \ref{HK props} \ref{HK linear} to \ref{HK Saks Henstock} for the HKP integral on integral currents of dimension $1$. 
Finally, we state three important convergence properties in the space of Henstock Kurzweil integrable functions:
\begin{Theorem}
  Let $(f_n)_n$ be a sequence of HK integrable functions on the interval $[\,a,b\,]$. Suppose that $f_n\to f$ pointwise almost everywhere. If any one of the following three conditions holds, then $f$ is HK integrable and $(HK)\int f =  \lim_n (HK)\int f_n$:
  \begin{enumerate}[label=(\roman*)]
    \item (Monotone Convergence Theorem) For almost all $x$, for all $n$, $f_n(x) \leq f_{n+1}(x)$ and  there holds $\sup_n (HK)\int f_n <+\infty$. 
    \item (Dominated Convergence Theorem) There exist HK integrable functions $g$ and $h$ such that for all $n$, $g\leq f_n\leq h$ almost everywhere. 
    \item (Controlled Convergence Theorem) $(f_n(x))_n$ is bounded for almost all $x\in [\,a,b\,]$ and for all $\epsilon>0$ there exists a positive gauge on $[\,a,b\,]$ such that for all $n$, for all $\delta$-fine tagged partition $\calP$ of $[\,a,b\,]$:
   \[
      \left \vert \sigma(\calP,f_n)-(HK)\int_a^b f_n \right \vert <\epsilon.
   \]
  \end{enumerate}
\end{Theorem}
In the latter case, the sequence $(f_n)_n$ is called HK \emph{equiintegrable}.
\begin{proof}
  The two first results can be proved using only the Saks Henstock Lemma and ``purely HK'' techniques, we will give such a proof for the Monotone Convergence Theorem of HKP Integration (see Theorem \ref{1d monotone convergence}). However, when possible, it is quicker to rely on Lebesgue integration results: the first statement follows from the Monotone Convergence Theorem of Lebesgue Integration. Indeed, if $f_1\leq f_n$ and both functions are HK integrable, then $f_n-f_1$ is nonnegative and HK integrable, thus Lebesgue integrable. To conclude, it suffices to apply Lebesgue's Monotone Convergence Theorem to the sequence $(f_n-f_1)$. Similarly, to prove the second result, consider the sequence $f_n-g$ and use the Lebesgue Dominated Convergence Theorem using $h-g$ as an upper bound.

The third statement has no equivalent in Lebesgue Integration, and relies strongly on the use of gauges. First redefine the $f_n$ and $f$ so that $f_n\to f$ everywhere and $(f_n(x))_n$ is bounded for all $x\in [\,a,b\,]$, this will not change the statement since the HK integral does not depend on the value of the function on a Lebesgue null set. Now, for $\epsilon>0$, choose $\delta$ as in the definition of the equiintegrability of the $f_n$. Let $\calP_1$ and $\calP_2$ be two $\delta$-fine tagged partitions of $[\,a,b\,]$, for all $n$, using the integrability of $f_n$ yields
\begin{align*}
  \left \vert \sigma(f,\calP_1) -\sigma(f,\calP_2) \right \vert \hspace{-2cm} \\
  & \leq\left \vert \sigma(f,\calP_1) -\sigma(f_n,\calP_1) \right \vert +  
   \left \vert \sigma(f_n,\calP_1) -\sigma(f_n,\calP_2) \right \vert\\
   &\hspace{6cm}+   \left \vert \sigma(f_n,\calP_2) -\sigma(f,\calP_2) \right \vert \\
&\leq \left \vert \sigma(f,\calP_1) -\sigma(f_n,\calP_1) \right \vert + 2\epsilon +\left \vert \sigma(f_n,\calP_2) -\sigma(f,\calP_2) \right \vert\\
&\leq  \sum_{(x,I)\in \calP_1} \vert f(x)-f_n(x)\vert \vert I\vert +2\epsilon +  \sum_{(x,I)\in \calP_2}  \vert f(x)-f_n(x)\vert \vert I\vert.
  \end{align*}
Since $\calP_1$ and $\calP_2$ are finite families, and $f_n$ converges to $f$ pointwise, for $n$ large enough we have
\[
   \left \vert \sigma(f,\calP_1)-\sigma(f,\calP_2)\right \vert < 3\epsilon
\]
and by the Cauchy criterion for HK integrability (see Proposition \ref{Cauchy HK}), $f$ is Henstock-Kurzweil integrable on $[\,a,b\,]$. To see that the integral of $f$ is the limit of the integrals of the $f_n$, consider $\delta$ adapted to $\epsilon$ for the integrability of the $f_n$ and for the integrability of $f$. Fix a $\delta$-fine tagged partition $\calP$. For $n$ large enough, there holds
   \begin{align*}
   \left \vert (HK)\int f - (HK)\int f_n \right \vert  \leq \hspace{-6cm} \\
 &\left \vert (HK)\int f -\sigma(f,\calP) \right \vert
+ \left \vert (HK)\int f_n -\sigma(f_n,\calP) \right \vert  +  \left \vert \sigma(f_n,\calP) -\sigma(f,\calP) \right \vert  \\
&\hspace{10cm}< 3\epsilon.
   \end{align*}
\end{proof}
\subsection{$AC_*$ functions and the Fundamental Theorem of Calculus}
The next paragraph follows closely the presentation of Sections 1.9 to 1.11 in T. De Pauw's survey \cite{DePauw2004autour}. We start with a first version of the Fundamental Theorem of Calculus for the Henstock-Kurzweil Integral.
\begin{Theorem}\label{HK First FTI}
  If $F$ is continuous on $[\,a,b\,]$ and differentiable at all but countably many points, then $F'$ is HK integrable on $[\,a,b\,]$ and $F$ is the indefinite integral of $F'$.
\end{Theorem}
\begin{proof}
Define $f$ to be equal to $F'$ wherever $F$ is differentiable and to $0$ elsewhere. Since $f$ is equal to $F'$ almost everywhere, $F'$ is HK integrable if and only $f$ is.
  
Fix $\epsilon>0$, let $y_1,y_2,\dots$ be the points at which $F$ is not differentiable. For $x\in [\,a,b\,]\backslash \{y_1,y_2,\dots\}$, using the differentiability of $F$ at $x$, choose a positive $\delta(x)$ such that for all $y\in [\, x-\delta(x),x+ \delta(x)\,]$
  \[
     \vert F(y) - F(x) - F'(x) (y-x)\vert <\dfrac{\epsilon}{2} \dfrac{\vert y-x\vert}{b-a}.
  \]
  For $j=1,2,\dots$, using the continuity of $F$ at $y_j$, choose $\delta(y_j)$ so that whenever $[\,c,d\,]$ is an interval in $[\,a,b\,]$ containing $y_j$ with $d-c <\delta(y_j)$, there holds
  \[
    \vert F(d)- F(c)\vert <\dfrac{\epsilon}{2^{j+2}}.
  \]
  Suppose $\calP =  (([\,c_k,d_k\,], x_k))_{k=1,\dots,p}$ is a $\delta$-fine tagged family in $[\,a,b\,]$. We can suppose up to reindexing $\calP$, that there exists $q\leq p$ such that for $k \in \{1,\dots, q\}$, $F$ is differentiable at $x_k$, whereas for $k\in \{q+1,\dots, p\}$, there exists $j$ such that $x_k= y_j$. Note also that a given $y_j$ corresponds to at most two different values of $k$ as no more than two nonoverlapping non trivial intervals can contain the same point.
  \begin{align*}
    \sum_{k=1}^p \left \vert F(d_k)-F(c_k) - f(x_k)(d_k-c_k)\right \vert \hspace{-6cm}\\
     & \leq \sum_{k=1}^q \left \vert F(d_k)-F(c_k) - F'(x_k)(d_k-c_k)\right \vert + \sum_{k=q+1}^p \left \vert F(d_k)-F(c_k) \right \vert\\
    &<\dfrac{\epsilon}{2} \dfrac{\sum_{k=1}^q (d_k-x + x-c_k)}{b-a} + 2\sum_{j=1}^\infty \dfrac{\epsilon}{2^{j+2}} \\
    &<\epsilon.
  \end{align*}
  Since $\epsilon>0$ is arbitrary, we apply this estimate to the case where $\calP$ is a partition to show that $f$ and $F'$ are HK integrable in $[\,a,b\,]$. As this is true for any tagged family, this shows that $F$ is the indefinite integral of $f$ and $F'$. 
\end{proof}
We now generalize this result to less regular functions $F$. This requires a new notion.
A function $F$ defined on $[\,a,b\,]$ is \emph{$AC_*$} if for every set $D\subseteq [\,a,b\,]$ of zero Lebesgue measure and every $\epsilon>0$, there exists a positive gauge $\delta$ on $D$ such that whenever $\calP$ is a $\delta$-fine family in $[\,a,b\,]$ tagged in $D$, there holds
\begin{equation}\label{ACstar eps small}
   \sum_{([\,c,d\,],x)\in \calP} \vert F(d)-F(c)\vert <\epsilon.
\end{equation}
In particular, an $AC_*$ function is continuous.
If $f$ is HK integrable, then its indefinite integral $F$ is $AC_*$, indeed if $D$ is a Lebesgue null set, we can consider the function $f_{D^c}\defeq f\ind_{D^c}$. 
As HK integration is insensitive to modifications on Lebesgue null sets, $F$ is also the primitive of $f_{D^c}$, so for $\epsilon>0$, we can apply the Saks-Henstock Lemma \ref{HK props}\ref{HK Saks Henstock} and find a gauge $\delta$ corresponding to $\epsilon/2$ on $[\,a,b\,]$. Considering the gauge $\delta_D = \delta\vert_D$ by the Saks-Henstock Lemma for any $\delta_D$ fine tagged family $\calP$ in $[\,a,b\,]$, since $f_D$ is equal to zero on $D$, \eqref{ACstar eps small} holds. The following converse statement holds:
\begin{Proposition}\label{ACstar HK}
  If $F$ is $AC_*$ and almost everywhere differentiable in $[\,a,b\,]$ then $F'$ is HK integrable and
  \[
    F(b)-F(a) = (HK)\int_a^b F'.
  \]
\end{Proposition}
\begin{proof}
  For $\epsilon>0$, define $\delta$ first on the set of differentiability points as in the previous proof and define $\delta$ on the null set of non differentiability points as a gauge adapted to $\epsilon/2$ in the definition of $AC_*$ functions.
\end{proof}
In order to get a general condition which ensures that a function is $AC_*$ and almost everywhere differentiable, we recall Stepanoff's Theorem. A function $F$ defined on an interval $I$ is \emph{pointwise Lipschitz} at the point $x\in I$ if
\[
    \Lip_x F \defeq \limsup_{y\to x y\in I} \dfrac{\vert F(y)-F(x)\vert }{\vert y-x\vert } <+\infty
\]
\begin{Theorem}[Stepanoff]\label{Stepanoff}
  If $F$ pointwise Lipschitz at all points of some set $E\subseteq I$, then $F$ is differentiable almost everywhere in $E$.
\end{Theorem}
For a proof of this result, see for instance \cite[Theorem 3.1.9]{FedererGMT}.
 \begin{Proposition}\label{ptwise Lip ACstar}
   A continous function $F$ which is pointwise Lipschitz at all but countably many points is $AC_*$.
 \end{Proposition}
 \begin{proof}
   Let $D$ be a null set in $[\,a,b\,]$ let $E \defeq \{y_1,y_2,\dots\}$ be the subset of $D$ at which $F$ is not pointwise Lipschitz. Fix $\epsilon>0$. Fix $k$ in $\{1,2,\dots\}$; as $F$ is continuous at $y_k$,  there exists $\delta_k>0$ such that whenever $x,y$ are points in $(y-\delta_k,y+\delta_k)\cap [\,a,b\,]$, 
\[
    F(y)-F(x) < 2^{-k-2}\epsilon.
\]
Now for $x\in D$, let
   \[
     \delta(x) \defeq \begin{cases}
       \dfrac{\epsilon}{2(b-a)(1+\Lip_x F)} \text{ if } x\in D\backslash E,\\
       \delta_k \text{ if } x= y_k \text{ for some } k \in \{1,2,\dots\}.
     \end{cases}
   \]
   The end of the proof is similar to that of Theorem \ref{HK First FTI}.
 \end{proof}
 Combining Proposition \ref{ptwise Lip ACstar}, Theorem \ref{Stepanoff} and Proposition \ref{ACstar HK} yields
\begin{Theorem}[Fundamental Theorem of Henstock Kurzweil Integration]\label{HK FTI}
  Let $F$ be a continuous function on $[\,a,b\,]$. Suppose that $F$ is pointwise Lipschitz at all but countably many points. Then $F$ is differentiable almost everywhere and $F'$ is Henstock Kurzweil integrable on $[\,a,b\,]$ with indefinite integral $F$.
\end{Theorem}
\subsection{An equivalent definition of the HK integral}
\begin{figure}
\begin{minipage}{.4\textwidth}
\includestandalone[width = \textwidth]{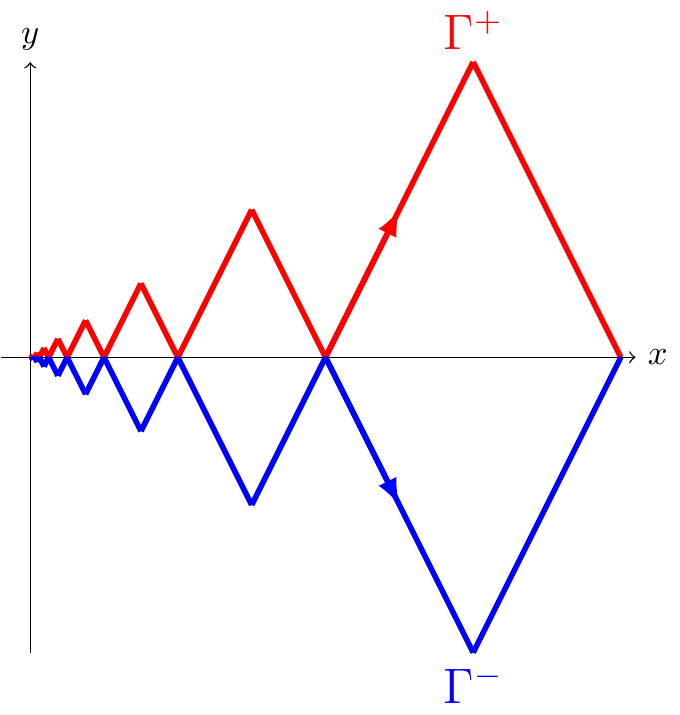}
\end{minipage}
\hspace{1cm}
\begin{minipage}{.4\textwidth}
\includestandalone[width = \textwidth]{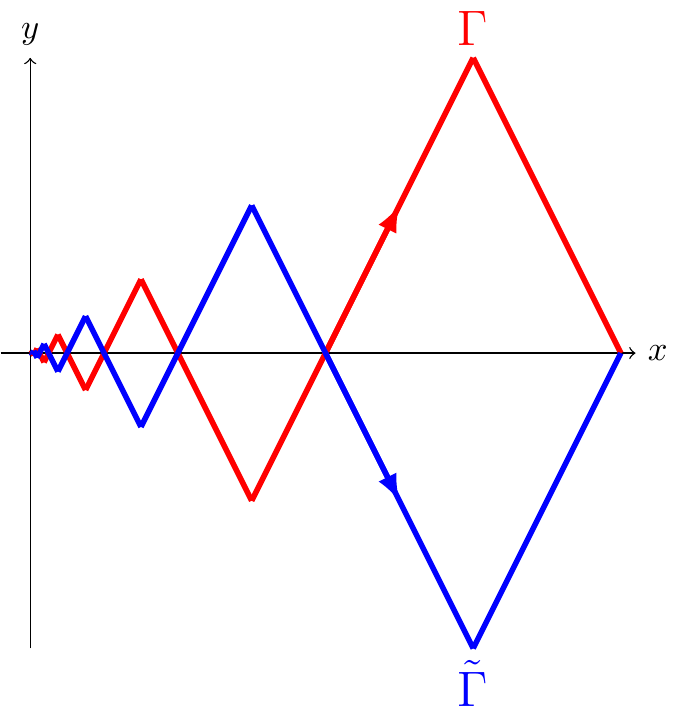}
\end{minipage}
\begin{centering}
\caption{$u$ is HK integrable on $\Gamma^+$ and $\Gamma^-$ but not on $\Gamma$ or $\tilde{\Gamma}$.}
\label{twocurvesfig}
\end{centering}
\end{figure}

\begin{Remark}[Extension to Lipschitz curves]
All the above properties of the Henstock Kurzweil Integral can be extended to the case where the interval $[\,a,b\,]$ is replaced by a simple Lipschitz curve $\Gamma\subseteq \Rn$ (closed or not). Indeed, one can consider an arc-length parameterization $\gamma$ of $\Gamma$ and work on $f\circ \gamma$. If $f$ is pointwise Lipschitz at $\gamma(x)$ along $\Gamma$, $f\circ \gamma$ is pointwise Lipschitz at $x$. The only thing that is not straighforward is relating differentiation in the ambient space $\Rn$ with differentiation along the curve. However, a Lipschitz curve has a tangent line at almost all points. In the next section, we consider countable sums of simple Lipschitz curves to develop Henstock-Kurzweil integration on integral currents of dimension $1$. The sum of curves can often be decomposed in several ways and Example \ref{example two curves} shows that the choice of the decomposition can have an effect on the integral, hence the need for a definition of integrability that does not depend on the decomposition.
\end{Remark}
\begin{Example}\label{example two curves}
  In $\R^2$, consider the curve $\Gamma^+$ corresponding to the graph in $(0,1\,]$ of the function
\[
     x\mapsto f(x)\defeq \dist (x, \{ t\in (0,1\,], 2t \sin (t^{-2}) - 2t^{-1} \cos(t^{-2})\}=0).
\]
 The curve  $\Gamma^+$ is a Lipschitz curve and has length $\sqrt 2$, orient $\Gamma^+$ towards the positive first coordinate. Let $ \Gamma^-$ be the the reflection of $\Gamma^+$ across the horizontal axis. The union of curves $\Gamma^+$ and $\Gamma^-$ can also be seen as the (closure of) the union of the graphs on $(0,1\,]$ of $x\mapsto \pm \sgn(x \sin (x^{-2}) - 2x^{-1} \cos(x^{-2})) f(x)$. Let $\Gamma$ and $ \tilde\Gamma $ be the corresponding curves. Let $u$ be the function defined in $\R^2$ by
  \[
    (x_1,x_2) \mapsto \begin{cases}
      2\sgn(x_2) \left(x_1 \sin (x_1^{-2}) - 2 x_1^{-1} \cos(x_1^{-2})\right) \hspace{-.2cm}&\text{ if } x_1>0, x_2\neq 0,\\
      \hspace{3cm}0 &\text{ otherwise.}
    \end{cases}
  \]
If $\gamma^+$, $\gamma^-$, $\gamma$ and $\tilde \gamma$ are respective arclength parametrizations of the curves above, the functions 
$u\circ \gamma^+$ and $u\circ \gamma^-$ are HK integrable on $[\,0, \sqrt 2\,]$ with respective indefinite integrals $x\mapsto \pm\sqrt{2}  x^2 \sin (x^{-2})$. However, the functions $u\circ \gamma$ and $u\circ \tilde\gamma$ are equal respectively to $\pm \vert (u\circ \gamma^+)'\vert$ which are not HK integrable. These curves are plotted in figure \ref{twocurvesfig}
\end{Example}
In order to generalize the Henstock Kurzweil integral to other settings, it is necessary to use more flexible tools. In particular we need to remove the dependency on the parameters and allow for families instead of partitions so that some ``small part'' of the domain can be left out. The precise meaning of a ``small part'' is a key point here.

This will be formalized in the next section, but first state an equivalent definition of HK integrability on an interval. In order to define what ``small'' is we will consider functions $F$ on the space of finite unions of disjoint intervals in $[\,a,b\,]$. Such a function is \emph{subadditive} if given two families, $\calU$ and $\calU'$, of closed intervals  of $[\,a,b\,]$ there holds
\[
    \left \vert F \left ([\calU] \cup [ \calU'] \right ) \right \vert \leq  \left \vert F \left ([\calU] \right ) \right \vert  +  \left \vert F \left ([\calU']\right ) \right \vert.
\]
$F$ is \emph{additive} if for $\calU$ and $\calU'$ as above with $\calL^1([\calU] \cap [\calU']) = 0$, there holds
\[
     F\left ([\calU]\cup [\calU'] \right ) = F ([\calU] ) + F ([\calU']).
\]
$F$ is \emph{continuous}  on the space of finite unions of intervals if given a sequence $\calU_j$ of families of intervals with $\# \calU_j <C$ and $\calL^1([\calU_j]) \to 0$, there holds $F([\calU_j]) \to 0$. In particular, it is equivalent to consider a continuous function $F:[\,a,b\,]\to \R$ and a continuous function on the space of finite unions of intervals of $[\,a,b\,]$, indeed a continuous function on a closed interval is uniformly continuous.

This definition seems impractical but we will see in the following section that it can be easily generalized to other supports and also to higher dimensions, in Section \ref{section pieces}. Indeed while intervals are not well suited to algebraic operations, they can be seen as currents of dimension $1$ in $\R$, using their canonical orientation and giving them multiplicity $1$. The following property is a reformulation of HK integrability in the language of Pfeffer integration (see Theorem 6.7.5 in \cite{Pfeffer1993riemannbook}).
\begin{Theorem}[Equivalent integrability condition]\label{HK equivalent def}
  A function $f$ defined almost everywhere on $[\,a,b\,]$ is Henstock-Kurzweil integrable if and only if there exists a nonnegative subadditive continuous function $G$ on the space of finite unions of intervals in $[\,a,b\,]$ and a real number $I$ with the property that for all $\epsilon>0$, there exists a gauge $\delta$ - not necessarily positive everywhere - and a positive number $\tau$ such that whenever $\calP$ is a $\delta$-fine tagged family in $[\,a,b\,]$ with
\[
    G([\,a,b\,] \backslash [\calP]) <\tau,
\]
there holds $\left \vert I-\sigma(f,\calP)\right \vert <\epsilon$.
\end{Theorem}
Before proving this, it makes sense to check that a tagged family satisfying the above contraints exists, it is a sort of generalization of Cousin's Lemma \ref{HK Cousin}, where we consider families instead of partitions.
\begin{Lemma}[Howard-Cousin Lemma]\label{HK Howard Cousin}
  Let $\delta$ be a gauge on $[\,a,b\,]$ (not necessarily positive). Let $G$ be a nonnegative, subadditive, continuous function on the space of finite union of closed intervals in $[\,a,b\,]$. For every $\tau>0$, there exists a $\delta$-fine tagged family $\calP$ in $[\,a,b\,]$ with
  \begin{equation}\label{HK G small remainder}
    G([\,a,b\,] \backslash [\calP]) <\tau.
  \end{equation}
\end{Lemma}
\begin{proof}
  We define a positive gauge $\delta^*$ on $[\,a,b\,]$ and use Cousin's Lemma to get a $\delta^*$-fine tagged partition $\calP^*$ of $[\,a,b\,]$ we then consider the subfamily $\calP$ of $\calP^*$ consisting of the pairs $(x,I)$ where $\delta(x)>0$. $\calP$ is clearly a $\delta$-fine tagged family, but it is necessary to check that \eqref{HK G small remainder} holds. This is where the choice of $\delta^*$ is critical. It relies 
heavily on the continuity of $G$:
    For every $\epsilon>0$, there exists $\eta>0$ such that whenever $I$ is a closed interval in $[\,a,b\,]$ with $b-a <\eta$, $G(I) < \epsilon$. To see this, suppose that the contrary holds: there exists a sequence $(I_j)_j$ of closed intervals in $[\,a,b\,]$ with $\calL^1(I_j)< 1/j$ and $G(I_j)\geq \epsilon$ for all $j$. This contradicts the hypothesis on $G$.

For $j=1,2,\dots,$ choose $\eta_j$ so that $G(I)<2^{-(j+1)}\tau$ whenever $I$ is an interval in $[\,a,b\,]$ with length less than $\eta_j$. Let $\delta^*(x_j) = \eta_j$ for $j=1,2,\dots$ and for $x\in [\,a,b\,]\backslash E$, let $\delta^*(x) \defeq \delta(x)$. $\delta^*$ is a positive gauge, so there exists a $\delta^*$-fine tagged partition $\calP^*$ of $[\,a,b\,]$. Clearly the tagged family $\calP \defeq \{(I,x)\in \calP^*, x\in [\,a,b\,]\backslash E\}$ is $\delta$-fine, and furthermore as $\calQ \defeq \calP^*\backslash \calP$ is finite and $G$ is subadditive, there holds
\[
    G([\,a,b\,]\backslash [\calP]) = G\left([\calP^*\backslash \calP]\right) \leq \sum_{(I,x)\in \calQ} G(I).
\]
Now for each $(I,x)\in \calQ$, there exists an index $j$ such that $x= x_j$ and there holds $\calL^1(I)<\delta^*(x_j)= \eta_j$. On the other hand, given an index $j$, there are at most two pairs of the form $(I,x_j)\in\calQ$.  Thus we have
\[
    G([\,a,b\,]\backslash [\calP]) < 2\sum_{j=1}^{\infty} 2^{-(j+1)}\tau \leq \tau.
\]
\end{proof}
We also need the following result on continuous, nonnegative, subadditive functions and gauges.
\begin{Lemma}\label{gaugeposgauge}
  Let $G$ be a continuous, nonnegative, subadditive function on the space of finite unions of subintervals of $[\,a,b\,]$ and $\delta$ be a gauge on $[\,a,b\,]$. Given $\tau>0$, there exists a positive gauge $\tilde{\delta}$ on $[\,a,b\,]$ which differs from $\delta$ only on the set $ E \defeq \{x\in [\,a,b\,], \delta(x) = 0\}$ and such that whenever $\tilde{\calP}$ is a $\tilde{\delta}$-fine tagged partition of $[\,a,b\,]$, the subfamily $\calP$ of $\tilde{\calP}$ consisting of the elements tagged in $[\,a,b\,]\backslash E$ satisfies
\[
    G([\,a,b\,]\backslash [\calP]) <\tau.
\] 
\end{Lemma}
\begin{proof}
  The set $E$ is countable, we can write $E= \{y_1,y_2,\dots\}$. For $j=1,2,\dots,$ by continuity of $G$, there exists $r_j>0$ such that if $I$ is an interval contained in $[\,a,b\,]$ and containing $y_j$, there holds $G(I)<2^{-j-1}\tau$. We thus define a positive gauge $\tilde{\delta}$ to be equal to $\delta$ in $[\,a,b\,]\backslash E$ and such that $\delta(y_j)= r_j$ for $j=1,2,\dots$. Consider a $\tilde{\delta}$-fine partition $\tilde{\calP}$ of $[\,a,b\,]$ and let $\calP$ be the tagged subfamily of $\tilde{\calP}$ consisting of the elements tagged outside in $[\,a,b\,]\backslash E$. As $\tilde{\calP}\backslash \calP$ is a finite family tagged in $E$ with no more than two elements tagged at the same point, we have
\[
    G([\,a,b\,]\backslash [\calP]) = G([\tilde{\calP}\backslash \calP]) < 2\sum_{j=1}^\infty 2^{-j-1} \tau \leq \tau.
\]
\end{proof}
\begin{proof}[Proof of Theorem \ref{HK equivalent def}]
First, we suppose that $f$ is HK integrable on $[\,a,b\,]$. Let $F$ be the indefinite integral of $f$. $F$ can be identified with an additive and continuous function on the space of finite unions of intervals in $[\,a,b\,]$. Let $G$ be the absolute value of $F$, it is subadditive, nonnegative and continuous. Fix $\epsilon>0$ and choose $\delta$ as in the definition of HK integrability for $f$. By the Saks-Henstock Lemma, if $\calP$ is a $\delta$-fine tagged family in $[\,a,b\,]$ there holds
\[
    \left \vert \sigma(f,\calP) -F([\calP]) \right\vert < 2\epsilon
\]
and by additivity of $F$
\[
    \left \vert F([\,a,b\,]) - F([\calP]) \right \vert \leq G([\,a,b\,] \backslash \calP),
\]
so if $\calP$ is such that
\[
   G([\,a,b\,]\backslash [\calP] ) <\epsilon,
\]
we have
\[
    \left \vert \sigma(f,\calP) -F([\calP]) \right\vert \leq\left \vert \sigma(f,\calP) -F([\calP]) \right\vert  + \left \vert F([\,a,b\,]) - F([\calP]) \right \vert \leq 3\epsilon.
\] 
This proves that $f$ satisfies the condition of the statement with $I \defeq (HK)\int f = F([\,a,b\,])$.

For the converse, we define $f$ everywhere in $[\,a,b\,]$. Fix $G$ and $I$ as in the statement and pick $\epsilon >0$. Chose a positive number $\tau$  and a gauge $\delta$ corresponding $\epsilon/2$ in the statement and let $E=\{y_1,y_2,\dots\}$ be the zero set of $\delta$. For $j=1,2,\dots,$ choose $r_j>0$ as in the proof of Lemma \ref{gaugeposgauge}. Taking $r_j$ smaller, we can also ask that $r_j \vert f(x_j)\vert <2^{-j-1} \epsilon$. Define $\tilde{\delta}$ as above. If $\tilde{\calP}$ is a $\tilde{\delta}$-fine tagged partition of $[\,a,b\,]$, we define the subfamily $\calP$ as before. Estimate \eqref{HK G small remainder} holds, and thus we also have
\[
     \vert I-\sigma(f,\calP) \vert < \epsilon/2.
\]
For the Riemann sum over $\tilde{\calP}$, we have:
\begin{align*}
    \vert I-\sigma(f,\tilde{\calP})\vert &\leq \vert I-\sigma(f,\calP) \vert + \sum_{(I,x)\in \tilde{\calP}\backslash \calP} \vert f(x)\vert \diam (I) \\
&<\epsilon/2 + 2\sum_{j=1}^\infty \vert f(y_j)\vert r_j \\
&< \epsilon.
\end{align*}
The function $f$ is thus HK integrable on $[\,a,b\,]$.
\end{proof}
\begin{Remark}
  The integral of Henstock and Kurzweil is \emph{not} equivalent to that of Pfeffer. The integrability condition in the above statement differs from that of Pfeffer in that the latter considers families made of regular sets of finite perimeter --- in $1$ dimension, these are finite unions of intervals. See Example 12.3.5 in \cite{Pfeffer1993riemannbook}.
\end{Remark}

\section{Integral currents of dimension $1$ and their pieces}\label{section pieces}
\subsection{Notations}
In the following, $f\vert_A$ denotes the restriction of the function $f$ to the set $A$, while $\mu \hel f$ and $\mu \hel A$ denote the multiplication of the (possibly vector valued) measure $\mu$ by the (scalar) function $f$ or the indicator function of $A$. $\spt \mu$ is the support of $\mu$.  In $\Rn$, with the usual Euclidean metric, we denote the norm of a vector $x$ by $\vert x\vert$ and the distance by $\dist(\cdot, \cdot)$. The usual scalar product of $x,y\in\Rn$ is $x\cdot y$, while the product of a vector $v$ with a covector $\eta$ is denoted $\langle \eta,v \rangle$.  $\Ball(x,r)$ and $\cBall(x,r)$ are respectively the open and closed balls of center $x\in \Rn$ and radius $r>0$. 

The Hausdorff measure of dimension $1$ is denoted by $\calH^1$. If $\mu$ is a scalar measure, $\set_1\mu$ denotes the points $x$ where $\mu$ has positive lower $1$-density, i.e. where
\[
    \Theta^{1*} (\mu,x) \defeq \limsup_{r\to 0} (2r)^{-1} \mu(\cBall(x,r)).
\]
A set $E\subseteq \Rn$ is \emph{$1$-rectifiable} if there exists a countable collection of Lipschitz curves $\gamma_j: \R\to \Rn$ such that $\calH^1(E\backslash\bigcup_j \gamma_j(\R)) =0$. A set is \emph{$0$-rectifiable} if it is countable.

An \emph{current} of dimension $m$ in $\Rn$ is a continuous functional on the space of smooth differential forms of degree $m$ with compact support: $\calD^m(\Rn)$. The space of such currents is denoted by $\calD_m(\Rn)$. In particular, a current of dimension $0$ is a distribution. The mass of a current $T\in \calD_m(\Rn)$ is the number
\[
   \mass(T)\defeq \sup \{ T(\omega), \omega\in \calD^m(\Rn), \Vert \omega\Vert \leq 1\} \in [0,+\infty].
\]
Currents of finite mass are representable by integration and $0$-currents with finite mass are measures. We denote by $\Vert T\Vert$ the carrying measure of a current. The \emph{boundary} of a current $T\in \calD_m(\Rn)$ for $m\geq 1$ is the current $\partial T \in \calD_{m-1}$ defined by $\partial T(\omega) = T(\dd \omega)$ for all $\omega\in \calD^{m-1}(\Rn)$. The \emph{flat norm} of a current $T\in \calD_{m}(\Rn)$ is the number
\[
   \flatn(T)\defeq \sup \{ \mass(Q)+\mass(R),\, T= \partial Q + R \} \in [0,+\infty].
\]

If $\gamma: [\,0,t_1\,] \to \Rn$ is a simple Lipschitz curve, we denote by $\lseg \gamma \rseg$ the current of dimension $1$ defined by
\[
    \lseg \gamma \rseg (\omega) = \int_0^{t_1} \langle \omega(\gamma(t)),\gamma'(t)\rangle \dd t.
\]
There holds $\mass(\lseg \gamma\rseg)  = \int_0^{t_1} \vert \gamma'(t) \vert \dd t$ and $\mass(\partial \lseg \gamma\rseg)$ is either $0$ or $2$ depending on whether $\gamma$ an open or a closed curve. The carrying measure of $\lseg \gamma \rseg$ is $\Vert \lseg \gamma \rseg \Vert = \calH^1 \hel \gamma([\,0,t_1\,])$.  We work mostly in \emph{integral currents of dimension $1$}, which include currents representing curves of finite length. A current $T\in \calD_1(\Rn)$ is integral ($T\in \Iirn$) if it can be written as a countable sum of simple Lipschitz curves $\lseg \gamma_j\rseg$, and there holds
\[
     \sum_j \mass(\lseg \gamma_j\rseg) = \mass (T) \text{ and } \sum_j \mass(\partial \lseg \gamma_j \rseg) = \mass (\partial T).
\]
In particular the density set $\set_1\Vert T\Vert$ is $1$-rectifiable. This characterization of integral currents is very specific to the one dimensional case.

 An integral current $T$ is \emph{decomposable} if there exists two non trivial integral currents $Q$ and $R$ with $Q+R=T$ and $\mass(T) = \mass(Q)+\mass(R)$, $\mass( \partial T )= \mass(\partial Q)+\mass (\partial R)$. If such a pair does not exist, $T$ is called indecomposable. 
A current $T\in \Iirn$ is indecomposable if and only if it is associated with an oriented simple Lipschitz curve with unit multiplicity.
\subsection{Pieces of a current}
\begin{figure}
  \begin{centering}
    \includestandalone[width =.4\textwidth]{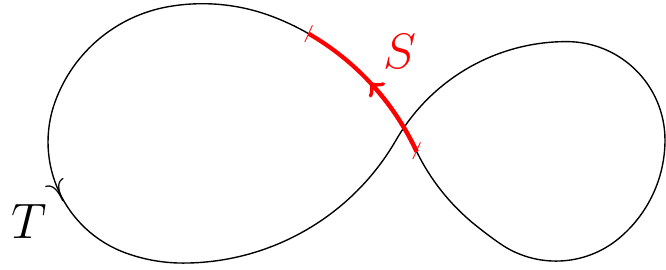}
    \caption{$S$ is not a piece of an indecomposable element of $T$.}
    \label{fig:notinadecomposition}
  \end{centering}
\end{figure}

Let $T$ be an integral current, an integral current $S$ is a \emph{piece of} $T$ if
\[
    \Vert S\Vert \leq \Vert T\Vert \text{ and } \Vert T-S\Vert \leq \Vert T\Vert.
\]
The notion of piece of a current differs from that of subcurrent defined in \cite{ThesisAJ,juliaARXIVstokes} for integral currents in any dimension where the condition is $\Vert S\Vert \perp \Vert T-S\Vert$. Subcurrents of $T$ are pieces of $T$, but the converse holds only if $T$ has multiplicity $1$ almost everywhere.
\begin{Example}
  Consider the current $T= 2 \lseg 0,2\rseg \in \mathbb{I}_1(\R^1)$, then
  \begin{itemize}
    \item The currents $\lseg 0,2\rseg$ and $2 \lseg 0,1\rseg$ are pieces of $T$,
    \item $3 \lseg 0,2\rseg$, $3^{-1} \lseg 0,2\rseg$ and $-\lseg 0,2\rseg$ are not pieces of $T$.
  \end{itemize}
\end{Example}
\begin{Proposition}
  An integral current $S$ is a piece of $T\in \mathbb{I}_1(\R^1)$ if and only if there exists a $\Vert T\Vert$ measurable function $g: \Rn\to [\,0,1\,]$ such that $S= T\hel g$.
\end{Proposition}
\begin{proof}
  Suppose $S= T\hel g$, then $\Vert S\Vert = \Vert T\Vert \hel g\leq \Vert T\Vert$ and $\Vert T-S\Vert = \Vert T\Vert \hel (1-g)\leq \Vert T\Vert$.

Conversely, suppose $S$ is a piece of $T$. Then $S$ is of the form $\calH^1\hel (\theta_S \ind_{M_S}) \wedge \vect S$ and $T = \calH^1\hel (\theta_T \ind_{M_T})\wedge \vect T$, where $\theta_S$ and $\theta_T$ are supposed non negative, respectively  $\calH^1\hel M_S$ and $\calH^1 \hel M_T$ almost everywhere. By the hypotheses on $S$ there holds
\begin{align*}
    \calH^1(M_S\backslash M_T)&= 0, \\
\theta_S &\leq \theta_T, \, \hspace{1cm} \calH^1\hel M_T \text{ almost everywhere},\\
 \vert \theta_T \vect T - \theta_S\vect S\vert  &\leq \theta_T, \, \hspace{1cm}\calH^1\hel M_T \text{ almost everywhere.}
\end{align*}
This in turn implies that $\vect T = \vect S$ at $\calH^1$ almost all points where $\theta_S$ is positive. Define the functions $g$ by
\[
   g(x) = \begin{cases}
 0 \text{ if } x\notin M_T \text{, or } \theta_T(x) = 0,\\
 \theta_S(x)/\theta_T(x) \text{ otherwise.}
\end{cases}
\]
Clearly $g(x) \in [\,0,1\,]$ for all $x\in \Rn$ and $S= T\hel g$.
\end{proof}
In particular elements of a decomposition of $T$ are pieces of $T$, however an indecomposable piece of $T$ may not be a piece of any element of decomposition of $T$ (see Figure \ref{fig:notinadecomposition}).
\subsection{Continuous function on the space of pieces of $T$}
Denote by $\subcspace_\leq (T)$ the collection of all pieces of $T$. 
\begin{Definition}  A function $F$ on $\subcspace_\leq(T)$ is \emph{continuous}, if given a sequence $(S_j)_j$ in $\subcspace_\leq(T)$ that converges to $0$ in the flat norm with $\sup_j \mass (\partial S_j) <+\infty$, we have $F(S_j) \to 0$. $F$ is \emph{additive} if whenever $S_1$ and $S_2$ are in $\subcspace_\leq (T)$ with $S_1+S_2 \in \subcspace_\leq(T)$ (which is equivalent to $\Vert S_1\Vert +\Vert S_2 \Vert \leq \Vert T\Vert$), there holds $F(S_1+S_2) = F(S_1)+F(S_2)$. $F$ is \emph{subadditive}, if instead for each $S_1$, $S_2$ as above, we have $F(S_1+S_2)\leq F(S_1)+F(S_2)$.
\end{Definition}
Continuous additive functions on $\subcspace_\leq(T)$ include the restriction of $1$-charges as defined in \cite{DPMoPf}: $1$-charges are defined on $\mathbb{N}_1(\Rn)$ and include in particular the representatives of continuous functions $f$ and continuous differential $1$-forms $\omega$ on $\spt T$ defined respectively as
\[
    \Theta_f : S\mapsto \partial S(f)
\]
and
\[
    \Lambda_\omega: S\mapsto \int \langle \omega, \vect S\rangle \dd \Vert S\Vert.
\]
Furthermore, the mass function $S\mapsto \mass(S)$ is continuous on $\subcspace_\leq(T)$:
%
%

%
%
\begin{Proposition}\label{1d mass continuous on pieces}
   For every $T\in \Iirn$, the function $S\mapsto \mass(S)$ is continuous and additive on $\subcspace_\leq(T)$.
 \end{Proposition}
 \begin{proof}
Additivity is clear. For the continuity, let $(S_j)_j$ be a sequence in $\subcspace_\leq(T)$ converging in the flat norm to $S\in \subcspace_\leq(T)$ with $\sup_j \mass( \partial S_j) <+\infty$. First notice that $\mass(S) \leq \liminf_j \mass(S_j)$ by lower semi-continuity of mass in the flat norm topology. So all we have to show is that $\liminf_j \mass(S_j)\geq \mass(S)$. In order to do this, for $\epsilon>0$ define a smooth $1$-form $\omega$ in $\Rn$ such that $\vert \omega(x)\vert \leq 1$ for all $x$ and $R(\omega) \geq \mass(R)-\epsilon$ for each $R\in \subcspace_\leq(T)$. 
%
%
Such a form exists, indeed by the definition of mass, there exists a smooth form $\omega$ such that  $\vert \omega(x)\vert \leq 1$ for all $x\in \Rn$ and $T(\omega) \geq \mass(T)-\epsilon$. Now, given $R\in \subcspace_\leq (T)$, there holds
\[
    R(\omega) = T(\omega) -(T-R)(\omega) \geq \mass(T) -\epsilon - \mass(T-R) \geq \mass(R) -\epsilon.
\] 
By definition of the flat convergence, $S_j(\tilde \omega)\to S(\tilde \omega)$ which implies that $\mass(S_j) \leq \mass(S) - \epsilon - \epsilon$, for all large enough $j$. Since $\epsilon$ is arbitary, $\mass(S_j)\to \mass(S)$.
 \end{proof}
As a consequence, to a $\Vert T\Vert$-Lebesgue integrable function $f$ in $\Rn$, one can associate the continuous additive function on $\subcspace_\leq(T)$: 
\[
   \tilde{\Lambda}_f : S\mapsto \int f \dd \Vert S\Vert.
\]
In the definition of $\Theta_f$, one can ask whether the continuity assumption of $f$ on $\spt T$ can be relaxed, for instance if $f$ is continuous on $\set_1 \Vert T\Vert$, is that sufficient for $\Theta_f$ to be continuous. Clearly, if $T$ is indecomposable, $\set_1\Vert T\Vert = \spt T$, but if one considers a current that has a countable decomposition, things are different:
\begin{Proposition}\label{example circles}
  There exists an integral current $T$ of dimension $1$ in $\R^2$ along with a bounded function $f$ continuous on $\set_1 \Vert T\Vert$, but not on $\spt T$ such that the function on $\subcspace_\leq(T)$ associated to the variation of $f$:
  \[ 
     \Theta_f : S\mapsto \partial S(f)
   \]
   is not continuous.
\end{Proposition}
\begin{proof}
  Consider a union of disjoint circles $\bigcup_{j=1}^\infty C_j$. Where for $j = 1,2,\dots$,  $C_j$ is centered at $(a_j,0) = (2^{-j},0)$ and has radius $r_j \defeq 3^{-(j+1)}$. Define the function $f$ piecewise on each $C_j$ so that $f= 1$ at the top (the point $(2^{-j},3^{-j-1})$) of each circle, and $f=-1$ at the bottom ($(2^{-j},-3^{-j-1})$) of each circle and $f$ is smooth. A good choice is $f(x_1,x_2) = r_j^{-1}y$ if $(x_1,x_2)\in C_j$. Let $\vect T$ be a field of tangent unit vectors to the circles, oriented positively and
\[
   T\defeq \left (\calH^2\hel \bigcup_j C_j \right) \wedge \vect T.
\]
Clearly $\spt T = \bigcup_j C_j \cup \{(0,0)\}$. Let us check that $\set_1\Vert T\Vert = \bigcup_j C_j$: for $r>0$ if $2^{1-j_0}\leq r\leq  2^{-j_0}$,
\[
    \Vert T\Vert(\Ball(0,r)) \leq \sum_{j\geq j_0} 2\pi r_j \leq  3^{j_0}\pi.
\]
Thus $\Theta^{1*}(\Vert T\Vert, 0)= 0$ and $0\notin \set_1 \Vert T\Vert$.
\begin{figure}
  \begin{centering}
    \includestandalone[width =.7\textwidth]{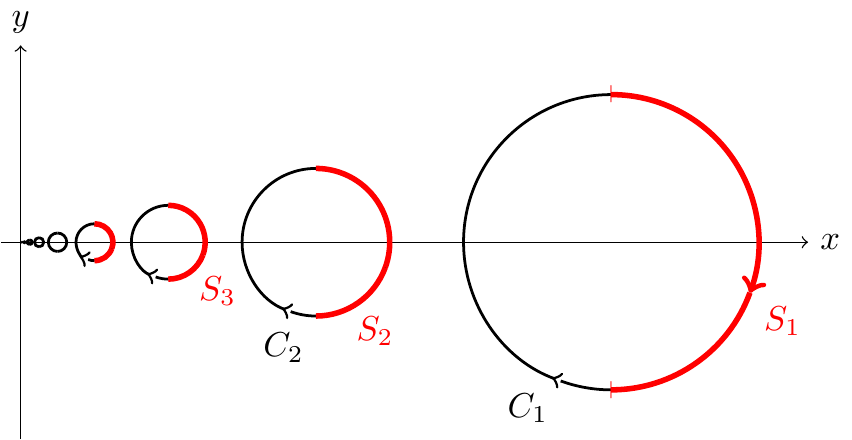}
    \caption{The current $T$ and the sequence $(S_j)_j$ of pieces}
    \label{fig:circles}
  \end{centering}
\end{figure}
  Consider the sequence of pieces $S_j\in \subcspace_\leq(T)$ corresponding to the half circles: $S_j = T\hel \{(x_1,x_2), 2^{-j}\leq x_1 \leq 2^{-j} + 3^{-j-1}\}$ (see figure \ref{fig:circles}). $S_j$ tends to $0$ in mass and for all $j$, $\mass(\partial S_j) = 2$. However, $\partial S_j(f) = 2 \nrightarrow 0$. Therefore the function $S\mapsto \partial S(f)$ is not continuous on $\subcspace_\leq(T)$.
\end{proof}
Note this never happens for an indecomposable current because of the following clear fact:
\begin{Fact*}
If $T\in \Iirn$ is indecomposable, then $\set_1\Vert T\Vert = \spt T$ and if $x\in \spt T\backslash \spt \partial T$, there holds $\Theta^1_*(\Vert T\Vert,x) \geq 1$.  
\end{Fact*}
\subsection{Derivation}
\begin{figure}
  \begin{centering}
    \includestandalone[width =.7\textwidth]{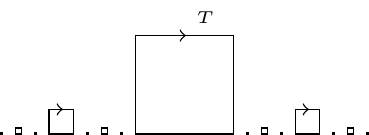}
    \caption{An integral current $T$ defined using the complementary intervals to a Cantor set.}
    \label{fig:cantor squares}
  \end{centering}
\end{figure}

 We use the terms derivation and derivate, following H. Federer \cite[Section~2.9]{FedererGMT}. For function on $\subcspace_\leq(T)$ there is a notion of derivation along $T$, similar to the differentiation of measures in Radon-Nikodym Theory:
\begin{Definition}
For $x$ in $\spt T$ and $\delta>0$, consider the subset $\subcspace_\leq(T,x,\delta)$ of $\subcspace_\leq(T)$ consisting of all pieces $S$ of $T$ such that 
  \begin{enumerate}
    \item $x\in \spt S$,
    \item $S$ is indecomposable,
    \item $\diam \spt S <\delta$.
  \end{enumerate}
If $\subcspace_\leq(T,x,\delta)$ is not empty for some positive $\delta$, the point $x$ is called \emph{good in $T$}.
In this case, we can define the \emph{upper and lower derivates} of $F$ along $T$ at $x$ respectively  as
\[
   \overline{\dcharge}_T F(x) \defeq \inf_{\delta>0} \sup_{\subcspace_\leq(T,x,\delta)}\dfrac{F(S)}{\mass(S)} \text{ and } \underline{\dcharge}_T F(x) \defeq \sup_{\delta>0} \inf_{\subcspace_\leq(T,x,\delta)}\dfrac{F(S)}{\mass(S)}.
\]
  $F$ is \emph{derivable along $T$ at} $x\in \set_1\Vert T\Vert$ if the upper and lower derivates of $F$ at $x$ along $T$ coincide, the corresponding \emph{derivate} is denoted $\dcharge_T F(x)$.

  A related notion we will use is that of almost derivability: a function $F$ on $\subcspace_\leq(T)$ is \emph{almost derivable} at $x\in \set_1 \Vert T\Vert$ if the upper and lower derivates of $F$ along $T$ at $x$ are finite.
\end{Definition}
We denote by $\textrm{Indec}(T)$ the set of points $x\in\Rn$ such that $\subcspace_\leq(T,x,\delta)$ is not empty for some $\delta>0$. $T$ has density at least $1/2$ at a point of $\textrm{Indec}(T)$, thus there holds 
\begin{align*}
    \textrm{Indec}(T) &\subseteq \set_1\Vert T\Vert\\
&\text{and}\\
\calH^1(\set_1\Vert T\Vert\backslash \textrm{Indec}(T)) &= 0.
\end{align*}
However, this last set can be large, as we show in the next proposition.
\begin{Proposition}
There exists an integral current $T$ of dimension $1$ in $\R^2$ such that $\set_1\Vert T\Vert\backslash \textrm{Indec}(T)$ is uncountable.
\end{Proposition}
\begin{proof}
  A way to define such a set is to consider a fat Cantor subset of $[\,0,1\,]$. For instance, one could let $C$ be the set obtained by removing iteratively the middle intervals of length $4^{-k}$ for $k=1,2,\dots$ from $[\,0,1\,]$. $C$ is a compact totally disconnected set with $\calL^1(C) =1/2>0$.

For each $k= 1,2,\dots$ there are $2^{k-1}$ segments of lenght $4^{-k}$ in the complement of $C$, denote them by $S_k^j$ for $j=1,2,\dots, 2^{k-1}.$ In $\R^2$ let $R_k^j$ be the rectangle $S_k^j\times[\,0,h_k\,]$ where the $h_k$ form a summable sequence of real numbers with $\sum_{k=1}^\infty 2^k h_k <+\infty$.

We can consider the current $T\in \mathbb{I}_1(\R^2)$ defined by
\[
   T\defeq \sum_{k,j} \lseg \bdry R_k^j \rseg.
\]
where the boundary curves of the squares are given a canonical orientation (see Figure \ref{fig:cantor squares}). $T$ is a cycle which has finite mass by the choice of $h_k$. Clearly $\spt T = C \cup \bigcup_{k,j} \bdry R_k^j \supseteq [\,0,1\,]$. The question is how to characterize $\set_1\Vert T\Vert$ and whether there exist points of $C\cap \set_1\Vert T\Vert$ such that there is no indecomposable piece $S$ of $T$ with $x\in \spt S$.
\begin{Claim}
Suppose $S$ is an indecomposable piece of $T$, then $S$ is a piece of $\lseg \bdry R_k^j\rseg $ for some $k\in \{1,2,\dots\}$ and $j\in\{1,\dots,2^{k-1}\}$.
\end{Claim}
\begin{innerproof}
  By contradiction, let $S\in \subcspace_\leq (T)$ be indecomposable and fix $x\in \spt S\cap \bdry R_k j$ and $x'\in \spt S\cap \bdry R_{k'}^{j'}$ with $(k,j)\neq (k',j')$. Without loss of generality (taking an indecomposable piece of $S$), we can suppose that $\partial S = \delta_{x'}-\delta_x$. We can also suppose that $x= (x_1,0)$ and $x'= (x'_1,0)$ with $x_1<x_1'$ and $x_1= \max(t \in S_k^j)$, $x'= \min(t, t\in S_{k'}^{j'}$. As $S$ is indecomposable and the differential form $(z_1,z_2)\mapsto \be_1^*$ is the differential of $(z_1,z_2)\mapsto z_1$, there holds
\[
    \int \langle \be_1^*,\vect S\rangle \dd \Vert S\Vert = x_1'-x_1.
\]
However, since $S$ is supported inside $[\,x_1,x'_1\,]\times \R$, and $\vect S= \vect T \in \{\be_1,-\be_1,\be_2,-\be_2\}$, $\Vert S\Vert$ almost everywhere, there holds
\begin{multline*}
    \int \langle \be_1^*,\vect S\rangle \dd \Vert S\Vert \leq \Vert T\Vert(([\,x_1,x'_1\,]\times \R)\cap\{(z_1,z_2), \vect T(z_1,z_2) = \be_1\}) \\
\leq \calL^1(C^c\cap [\,x_1,x'_1\,])< x'_1-x_1,
\end{multline*}
where we used the fact that $C\cap[\,x_1,x'_1\,]$ contains a fat Cantor subset of $C$, which has positive Lebesgue measure. This is a contradiction.
\end{innerproof}
The above claim implies that for all $x\in C\backslash \bigcup_{k,j} \cl(S_k^j)$, $x$ is not in the support of any indecomposable piece of $T$. There remains to prove that $C\cap \set_1\Vert T\Vert \backslash  \bigcup_{k,j} \cl(S_k^j)$ is uncountable.
For $x\in (0,1)$, $\Theta^{1*}(\Vert T\Vert,x) \geq \Theta^{*1}(\calH^1\hel E^c,x) = 1- \Theta^{1}_*(\calH^1\hel C,x)$, so we only need to prove that $C\cap \{x, \Theta^1_*(\calH^1\hel C,x)<1\}$ is uncountable.

In \cite[Theorem 1]{Buczolich2003category}, Buczolich proved that the set of points of a nowhere dense perfect set $P\subseteq \R$ where $P$ has lower density larger than $\gamma$ for any $\gamma>0.5$ is always of first category in $P$. This implies that the set of points of density less than $1$ is of second category in $P$, which in turn implies that it is uncountable ($P$ is a Baire space with the topology inherited from $\R$, see for instance \cite[Chapter 9]{Oxtoby1980MnC}). Note that there are more precise ways to characterize the points of a Cantor set with given densities, see for instance the paper by Besicovitch \cite{Besicovitch1956density}.
\end{proof}
\begin{Example}
   Let $\Lambda_f$ be the function on $\subcspace_\leq(T)$ be associated to a Lebesgue $\Vert T\Vert$ integrable function $f$ defined almost everywhere on $\set_1\Vert T\Vert$ by
\[
    \Lambda_f : S\mapsto \int f\dd \Vert S\Vert.
\]
If $f$ is continuous at $x\in \set_1\Vert T\Vert$ and $x$ is good in $T$, then $\Lambda_f$ is derivable at $x$ along $T$ with derivate $\dcharge_T F(x) = f(x)$. 
\end{Example}
For a good point $x\in \spt T$, $\epsilon>0$, choose $\delta >0$ such that $\vert f(y)-f(x) \vert < \epsilon$ for all $y\in \Ball(x,\delta)$. For $S\in \subcspace_\leq(T,x,\delta)$
\[
   \vert \Lambda_f(S)-f(x)\mass(S) \vert  \leq \int \vert f(y)-f(x)\vert \dd \Vert S\Vert(y) \leq \epsilon \mass(S).
\]
Letting $\epsilon$ go to zero, we can conclude.
\begin{Question}
  If $F$ is a continuous function defined on $\subcspace_\leq(T)$, are the extended real valued functions $\underline\dcharge_TF$, $\overline\dcharge_TF$ and $\dcharge_TF$ are $\Vert T\Vert$ measurable? Borel measurable?.
\end{Question}
For Henstock-Kurzweil Integration in $1$ dimension and for Pfeffer Integration on sets of finite perimeter, such results rely on the Vitali covering theorem and a derivation operation. A ``covering'' theorem using pieces of $T$ would be useful. An alternative would be to study a suitable decomposition of $T$, but this approach is made difficult by the fact that there can be pieces of $T$ which do not belong to any decomposition of $T$, as demonstrated in Figure \ref{fig:notinadecomposition}.

%
\begin{Definition}
  Let $T$ be an integral current of dimension $1$ in $\Rn$ and let $u$ be a function defined on $\set_1 \Vert T\Vert$. Fix a good point $x\in \set_1 \Vert T\Vert$. The function $u$ is \emph{differentiable along $T$} at $x$ if there exists a linear form $\dD u (x)$ on $\Rn$ such that for all $\epsilon>0$, there exists $\delta>0$ such that whenever $y\in \set_1 \Vert T\Vert \cap \Ball(x,\delta)$ and there is an $S\in \subcspace_\leq(T,x,3\delta)$ with $y\in \set_1 \Vert S \Vert$, there holds
\[
    \vert u(y)-u(x) -\dD u(x)\cdot (y-x) \vert \leq \epsilon \vert y-x\vert.
\]
\end{Definition}
Note that if $u$ is differentiable in $\Rn$ or differentiable on $\spt T$ in the sense of Whitney \cite{Whitney34} then $u$ is differentiable along $T$ with the same differential.
\begin{Theorem}\label{1d thm derivability}
   Suppose that $u$ is a continuous function on $\spt T$ for some $T\in \Iirn$. Fix $x\in \set_1 \Vert T\Vert$ such that $\subcspace_\leq(T,x,\delta)\neq \emptyset$ for some $\delta>0$, then the following three statements hold
   \begin{enumerate}[label =(\roman*)]
   \item \label{1d Du zero} If $u$ has pointwise Lipschitz constant $\Lip_x u= 0$ at $x$, then $\Theta_u$ is derivable at $x$ along $T$ and $\dcharge_T \Theta_u(x) = 0$.
   \item \label{1d u pointwise lip} If $u$ is pointwise Lipschitz at $x$, then $\Theta_u$ is almost derivable at $x$ with $- \Lip_x u \leq \underline \dcharge_T \Theta_u\leq \overline \dcharge_T \Theta_u \leq Lip_x u$.
     \item \label{1d u derivable} If $\vect T$ has a $\Vert T\Vert$ approximately continuous representative at $x$ (which we still denote by $\vect T$), $\Vert T\Vert$ has finite upper density at $x$ and $u$ is differentiable at $x$ along $T$, then $\Theta_u$ is derivable at $x$ along $T$, with $\dcharge_T \Theta_u(x) = \langle \dD u(x),\vect T(x)\rangle$.
   \end{enumerate}
 \end{Theorem}
 \begin{Remark}
   The assumption that $\vect T$ has a $\Vert T\Vert$ approximately continuous representative at $x$ is satisfied for $\Vert T \Vert$ almost all $x$. (See Claim \ref{vect T approx cont} in the proof of Proposition \ref{1d prop almost diff ACstar}.)
 \end{Remark}
 \begin{proof}
   Let us start with \ref{1d Du zero} and \ref{1d u pointwise lip}. For $\epsilon>0$, there exists $\delta$ such that whenever $y\in \spt T$ with $\vert y-x\vert<\delta$,
  \[
      \vert u(y)-u(x) \vert <(M+\epsilon) \vert y-x\vert,
    \]
    with $M\defeq \Lip_x u$. Given an indecomposable $S\in \subcspace_\leq(T)$, with $x\in \spt S$ and $\diam \spt S <\delta$, $S$ is of the form $\gamma_\# \lseg 0,\mass(S)\rseg$ with $\gamma(0) = y_-$ and $\gamma(\mass(S)) = y_+$. Since $\vert y_+ - x\vert + \vert x-y_-\vert \leq \mass(S)$, we get
      \begin{align*}
        \vert \Theta_u(S)\vert &= \vert u(y_+)-u(y_-)\vert \\ &\leq \vert u(y_+)-u(x)\vert  + \vert u(x)-u(y_-)\vert\\ &\leq  (M+\epsilon) \mass(S).
      \end{align*}
      As $\epsilon$ is arbitrary, this is enough to prove \ref{1d u pointwise lip}, and \ref{1d Du zero} where we have $M=0$.
   We turn to \ref{1d u derivable}. 

If $\dD u (x) = 0$, refer to \ref{1d Du zero}, thus we can suppose $\dD u (x) \neq 0$. Fix $\epsilon>0$. There exists $\delta_1>0$ such that for any $r\in (0,\delta_1)$,
\begin{equation}\label{eq: 1d T density}
     \dfrac{\Vert T\Vert(\cBall(x,r))}{2r} \leq 2\theta,
\end{equation}
with $\theta \defeq \Theta^{1*}(\Vert T\Vert,x)\in (0,+\infty)$. Replace $\vect T$ with its $\Vert T \Vert$ approximately continuous representative at $x$. Denote by $E_{x,\epsilon}$ the set
\[
   E_{x,\epsilon} \defeq \set_1\Vert T \Vert \cap \left \{ y, \vert \vect T(y)-\vect T(x)\vert > \dfrac{\epsilon}{2\vert \dD u (x) \vert} \right \}.
\]
 There exists $\delta_2>0$ which we can suppose less or equal to $\delta_1$ such that whenever $r\in(0,\delta_2)$,
\begin{equation}\label{eq: 1d T tangent}
    \dfrac{\Vert T\Vert (\cBall(x,r) \cap E_{x,\epsilon})}{\Vert T\Vert(\cBall(x,r))}
 <\dfrac{\epsilon }{4\theta \vert \dD u(x)\vert}.
\end{equation}
For $S\in \subcspace_\leq(T,x,\delta_2)$, the field $\vect S$ is equal $\Vert S\Vert$ almost everywhere to $\vect T$ and if $S$ represents a curve joining $x$ and $y$, with $\partial S = \delta_y-\delta_x$. As for $j=1,\dots,n$ the $1$ form $z\mapsto \be_j^*$ is the differential of the $0$-form $z\mapsto z_j$, We can write:
\begin{multline*}
    y-x = (y_1-x_1)\be_1 +\dots + (y_n-x_n)\be_n = \sum_{j=1}^n \partial S(z\mapsto z_j)\be_j\\
= \sum_{j=1}^n S(z\mapsto \be_j^*)\be_j = \sum_{j=1}^n \int \langle \be_j^*, \vect T\rangle \dd \Vert S\Vert \be_j = \int \vect T \dd \Vert S\Vert.
\end{multline*}
The same identity with opposite sign is true if $\partial S = \delta_x-\delta_y$ instead. Denote by $d_S$ the diameter of $\spt S$. By \eqref{eq: 1d T density} and \eqref{eq: 1d T tangent},
\begin{align}\label{eq: 1d S tangent}
   \vert y-x - \mass(S) \vect T(x)\vert  &\leq \int \vert \vect T(x')-\vect T(x) \vert \dd \Vert S\Vert(x')\\
&\leq 2 \Vert S\Vert \left (E_{x,\epsilon} \cap \cBall(x,d_S) \right ) + \dfrac{\epsilon }{2\vert \dD u(x)\vert}\mass(S)\nonumber\\
&\leq \dfrac{\epsilon \Vert T\Vert (\cBall(x,d_s)) }{2\theta \vert \dD u(x)\vert} + \dfrac{\epsilon }{2\vert \dD u(x)\vert}\mass(S)\nonumber\\
&\leq \dfrac{ 2\epsilon \theta d_S}{\theta \vert \dD u(x)\vert} + \dfrac{\epsilon }{2\vert \dD u(x)\vert}\mass(S)\nonumber\\
&\leq \dfrac{5 \epsilon}{2\vert \dD u(x)\vert} \mass(S),\nonumber
\end{align}
where in the second inequality, we used the fact that  $\vert \vect T(x')-\vect T(x)\vert \leq 2$ for $\Vert T \Vert $ almost all $x'$, in particular in the exceptionnal set $E_{x,\epsilon}$, in the third inequality we used the fact that $\Vert S\Vert \leq \Vert T\Vert$ and \eqref{eq: 1d T tangent} and in the last inequality, we used the fact that since $S$ is indecomposable, $d_S \leq \mass(S)$.

By differentiability of $u$ along $T$ at $x$, there exists $\delta_3>0$, such that for $y\in \Ball(0,\delta_3)\cap \set_1 \Vert T \Vert$ such that there exists $S\in \subcspace_\leq(T,x,\delta_3)$ with $y\in \spt S$,
\[
    \vert u(y)-u(x) - \langle  \dD u (x),  y-x\rangle \vert < \epsilon \vert y-x\vert.
\]
Let $\delta \defeq \min\{ \delta_1, \delta_2,\delta_3\}$ and choose $S \in \subcspace_\leq(T,x,\delta)$. We can write $S$ as $S^++ S^-$ where $S^+$ and $S^-$ are indecomposable, $\partial S^+= \delta_ {y^+} -\delta_x$ and $\partial S^- = \delta_x- \delta_{y^-}$, with $\mass(S) = \mass(S^+)+ \mass(S^-)$ and we have
\[
    \Theta_u(S)= \Theta_u(S^+) +\Theta_u(S^-) = u(y^+) - u(x) + u(x)-u (y^-).    
\]
Thus we can write
\begin{align*}
    \vert \Theta_u(S) -\langle \dD u (x) ,\vect T(x) \rangle \mass(S) \vert \hspace{-1cm}\\
&\leq \vert u(y^+)-u(x) -\langle \dD u (x) ,\vect T(x) \rangle \mass(S^+)\vert \\
&+ \vert u(x)-u(y^-) -\langle \dD u (x) ,\vect T(x) \rangle \mass(S^-)\vert
\end{align*}
  and study only the first term of the right hand side. We have
\begin{align*}
\vert u(y^+)-u(x) -\langle \dD u (x) ,\vect T(x) \rangle \mass(S^+)\vert \hspace{-5cm} \\
&\leq 
    \vert u(y^+)-u(x) -\langle\dD u(x), y-x\rangle\vert  \\
&\hspace{2cm}+ \vert \langle\dD u(x) ,y^+-x\rangle-\langle \dD u (x) ,\vect T(x) \rangle \mass(S^+)\vert\\
&\leq \epsilon \vert y^+-x\vert + \vert \dD u(x)\vert \vert y^+-x - \mass(S^+) \vect T(x)\vert \\
&\leq 4 \epsilon \mass(S^+),
  \end{align*}
by \eqref{eq: 1d S tangent} applied to $S^+$. Doing the same with $S^-$ and summing concludes the proof: there exists $\delta>0$ such that for all $S\in \subcspace_\leq(T,x,\delta)$,
\[
   \vert \Theta_u (S) -\mass(S) \langle \dD u (x),\vect T(x)\rangle \vert \leq \epsilon \mass(S)
\]
and $\Theta_u$ is thus differentiable along $T$ at $x$.
\end{proof}
If one assumes only approximate continuity of the tangent - as we just did - the assumption that the currents $S$ used in the derivation are indecomposable is necessary:
\begin{Example}
  Consider the function $h: (x,y)\mapsto y$ and the current $T$ associated to an infinite staircase with steps indexed by $j$, with height ($y$ length) $3^{-j}$ and length ($x$-length) $2^{-j}$ symmetric in the $x$ direction, converging at $(0,0)$ (see Figure \ref{fig:zigzag}). If one considers a sequence of subcurrents $S_j$ composed of a very small ``interval'' (length $4^{-j}$) around $0$ and a vertical part of the step, there holds
\[
     \Theta_h (S_j) = 3^{-j}C + o(3^{-j}).
\]
Thus $\lim_j \Theta_h(S_j)/\mass(S_j) = C>0$.
\begin{figure}
  \begin{centering}
    \includestandalone[width =.7\textwidth]{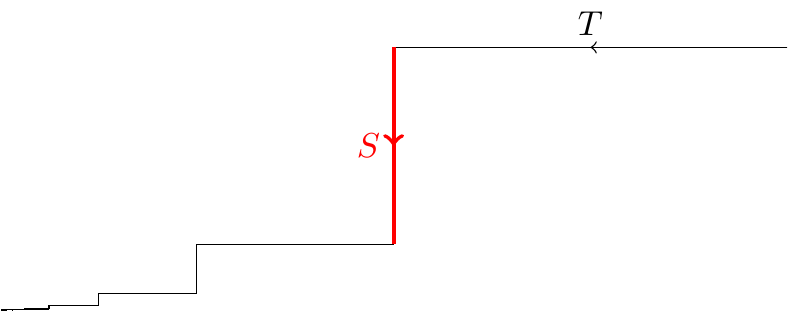}
    \caption{The piece $S$ is not suitable for a differentiation basis of $T$ at $0$.}
    \label{fig:zigzag}
  \end{centering}
\end{figure}
However if one considers a sequence of indecomposable currents $R_j$ touching $0$ with mass going to $0$, we will get by the above Theorem $\Theta_h(R_j)/\mass(R_j) \to 0$.
\end{Example}
An alternative restriction would be to bound the \textit{regularity} of the pieces. This is actually how we proceed in higher dimension in \cite{ThesisAJ, juliaARXIVstokes} as indecomposability is not a practical tool for currents of higher dimension.
\section{Integration}\label{sec 1d integration}
We first need an analogue to Cousin's Lemma in order to decompose a current of dimension $1$ into small pieces.
\subsection{Howard Cousin Lemma in dimension $1$}
Given a current $T\in \Iirn$ and a gauge on $\set_1 \Vert T \Vert$, a \emph{tagged family in $T$} is a finite collection $\calP$ of pairs $(S_j,x_j)$ for $j=1,\dots,p$,  where
\begin{align*}
   &S_j \in \subcspace_\leq(T),\text{ is indecomposable},\\
   &x_j \in \set_1 \Vert T\Vert \cap \spt S_j \\
    &\text{and}\\
   &\sum_{j=1}^p \Vert S_j\Vert \leq \Vert T\Vert.
\end{align*}
 If $T$ has multiplicity $1$ almost everywhere, the last condition prevents the pieces from overlapping. Such a tagged family is \emph{subordinate to a decomposition}  $T_1,T_2,\dots$ of $T$ if there exists a partition of $\calP$ indexed by $k$ into families $\calP_k$ each in the respective $T_k$.

A \emph{gauge} on a set $E$ is a nonnegative function $\delta$ such that $\{ x \in E,\delta(x)=0\}$ is countable. If $\delta$ is a gauge on a set $E\subseteq \set_1 \Vert T\Vert$, a \emph{ $\delta$-fine} tagged family in $T$ is a tagged family as above satisfying
\[
   \forall (S,x) \in \calP,\,  x\in E \text{ and } \diam \spt S <\delta(x).
 \]
Furthermore, given a nonnegative subadditive function $G$ on $\subcspace_\leq(T)$, and a positive real number $\tau$, a tagged family $\calP$ is \emph{$(G,\tau)$ full} if $G(T-[\calP]) <\tau$. 
\begin{Lemma}[Howard-Cousin Lemma]\label{1dcurrent HC}
  Let $T$ be an integral current of dimension $1$ in $\Rn$. Let $F$ be a subadditive continuous function on $\subcspace_\leq(T)$. Given $\epsilon >0$  and $\delta$ a gauge on $\set_1 \Vert T\Vert$, for any decomposition $T_1,T_2,\dots$, there exists a $(\vert F\vert,\epsilon)$ full, $\delta$-fine tagged family subordinate to this decomposition.
\end{Lemma}
\begin{proof}
Fix a decomposition of $T$. For each $k$ choose $\gamma_k:[\,0,\mass(T_k)\, ] \to \Rn$ to parameterize $T_k$ by arc-length, so that $T_k = \gamma_{k\#} \lseg 0,\mass(T_k)\rseg$. Let $\delta_k\defeq \delta\circ \gamma_k$, it is a gauge on $I_k\defeq [\,0,\mass(T_k)\,]$.

Since $T$ is integral, there exists $k_0$ such that for all $k> k_0$, $T_k$ is a cycle. Also $\mass(T_k)\to 0$ as $k \to \infty$. Since $F$ is continuous and subadditive, there exists $k_\epsilon$ such that
\[
   \left \vert  F\left (T-\sum_{k= 1}^{k_\epsilon} T_k\right )\right \vert <\dfrac{\epsilon}{2}.
\]
For $k=1,2,\dots,k_0$, consider the interval $I_k[\,0,\mass(T_k)\,]$, along with the gauge $\delta_k$ and the continuous additive function $\gamma_k^{\#}F$ on $\subcspace_\leq(\lseg I_k\rseg)$ defined by $\gamma_k^{\#}F(\lseg a,b\rseg) = F(\gamma_{k\#} (\lseg a,b\rseg)$, for $0\leq a< b\leq \mass(T)$. Note that it is enough to define $\gamma_k^\# F$ on indecomposable pieces of $\lseg I_k \rseg$ as all pieces  are in this case a finite sum of disjoint indecomposable pieces. Apply Lemma \ref{HK Howard Cousin} to $I_k$, $\delta\circ \gamma_k$, $\vert \gamma_k^\# F\vert$ and $\epsilon/(2 k_0)$ to get a $\delta\circ \gamma_k$ fine $(\gamma_k^\#, \epsilon/(2k_0))$ full tagged family $\calP_k$ in $I_k$.

The collection $\gamma_\# \calP_k$ defined by $\{ (\gamma_{k\#} S,\gamma_k(x)), (S,x)\in \calP_k\}$ is a $\delta$ fine tagged family in $T_k$ (as $\gamma_k$ has Lipschitz constant $1$), which satisfies
\[
   \vert F(T_k- [\gamma_\#\calP_k]) \vert = (\gamma_k^\# F)(\lseg I_k \rseg - [\calP_k]) \vert < \dfrac{\epsilon}{2k_0}.
 \]
 Summing this inequality over $k=1,2,\dots,k_0$ yields
 \[
   \left \vert F \left (\sum_{k=1}^{k_0} T_k - \left [\bigcup_{k=1}^{k_0} \gamma_{k\#}\calP_k\right]\right )\right \vert < \dfrac{\epsilon}{2}.
 \]
And the collection $\calP \defeq \bigcup_{k=1}^{k_0} \gamma_{k\#}\calP_k$ is therefore a tagged family in $T$ which is $\delta$ fine and $(F,\epsilon)$ full.
\end{proof}
\subsection{$AC_*$ functions on $\subcspace_\leq(T)$}
A function $F$ on $\subcspace_\leq(T)$ is \emph{$AC_*$} if given a $\Vert T\Vert$ null set $E\subset \set_1\Vert T\Vert$, for every $\epsilon>0$, there exists a gauge $\delta$ on $E$ with
\[  
    \vert F([\calP])\vert <\epsilon,
\]
whenever $\calP$ is a $\delta$-fine tagged family in $T$. We say that a tagged family is \emph{anchored} in a set $E$ if for all $(S,x)$ in this tagged family, $x\in E$. As the gage $\delta$ is defined only on $E$, here $\calP$ is automatically anchored in $E$.
The next two propositions are adapted from of \cite[Theorems 3.6.6. and 3.6.7]{Pfeffer2001book}.
\begin{Proposition}
If $F$ is a continuous additive function on $\subcspace_\leq(T)$ which is $AC_*$ and such that $\underline{\dcharge}_T F(x)\geq 0$ almost everywhere, then $F$ is nonnegative, i.e. for all $S\in \subcspace_\leq(T)$, $F(S)\geq 0$.
\end{Proposition}
\begin{proof}
It is enough to prove that $F(T)\geq 0$, indeed if $T'$ is in $\subcspace_\leq(T)$ the restriction of $F$ to $\subcspace_\leq(T')$ satisfies the hypothesis of the proposition. Let $N$ be the set of points $x$ such that $\underline \dcharge_T F(x) <0$. For $\epsilon>0$, there exists a gauge $\delta_N$ on $N$ such that $\vert F([\calP])\vert <\epsilon$ whenever $\calP$ is a $\delta_N$ fine tagged family anchored in $N$. For each $x$ at which $\underline{\dcharge}_TF(x)\geq 0$, there exists $\Delta_x$ such that for all $S\in \subcspace_\leq(T,x,\Delta_x)$, $F(S)\geq -\epsilon \mass(S)/\mass(T)$. Define a gauge $\delta$ on $\set_1\Vert T\Vert$ by letting
\[
   \delta(x) = \begin{cases}
\delta_N(x) \text{ if } x\in N,\\
\Delta_x \text{ otherwise}.
\end{cases}
\]
Using Lemma \ref{1dcurrent HC}, find a $\delta$ fine tagged family $\calP$ in $T$ with $\vert F(T-[\calP]) \vert <\epsilon$. Let $\calP_N$ be the subfamily of $\calP$ consisting of all the elements anchored in $N$. Denoting $\calP^*$ the complement of $\calP$ yields:
\[
    F(T) \geq F([\calP]) - F(T-[\calP] )\geq F([\calP^*]) + F([\calP_N]) - \epsilon \geq -3\epsilon.
\]
Since $\epsilon$ is arbitrary, $F(T) \geq 0$.
\end{proof}
\begin{Proposition}\label{1d prop almost diff ACstar}
   If a continuous additive function $F$ is almost derivable everywhere in $\set_1\Vert T\Vert$  except in a countable set $E_T$, then $F$ is $AC_*$.
\end{Proposition}
\begin{proof}
Let $N$ be a $\Vert T\Vert$ null set. For $\epsilon>0$, and $k=1,2,\dots$, let $U_k$ be a neighbourhood of $N$ with $\Vert T\Vert (U_k) <2^{-k}\epsilon/k$. For $x\in N\backslash E_T$, choose a positive integer $k_x$ and a positive $\Delta_x$ such that $\Ball(x,\Delta_x) \subseteq U_{k_x}$ and for all $S\in \subcspace_\leq(T,x,\Delta_x)$, $\vert F(S)\vert\leq k_x\mass(S)$. $k_x$ and $\Delta_x$ exist by almost derivability of $F$ at $x$. Define a gauge $\delta$ on $N$ by
\[
   \delta(x) = \begin{cases}
0 \text{ if } x\in E_T,\\
\Delta_x \text{ if } x\in N\backslash E_T.
\end{cases}
\]
Given a $\delta$ fine tagged family $\calP$ anchored in $N$, partition $\calP$ into families $\calP_k$ for $k=1,2,\dots$ such that $(S,x)\in \calP_k$ if and only if $k_x = k$ all but finitely many of these families are empty, there holds
\begin{equation*}
 \vert F([\calP]) \vert \leq \sum_{k=1}^\infty \sum_{(S,x)\in \calP_k} \vert F(S)\vert \leq \sum_{k=1}^\infty k \sum_{(S,x)\in\calP_k}\mass(S) \leq \sum_{k=1}^\infty k\Vert T \Vert(U_k) < \epsilon.
\end{equation*}
\end{proof}
%
\subsection{The HKP Integral on integral currents of dimension $1$.}
\begin{Definition}\label{Pfeffer1 integrable}
  A function $f$ defined $\Vert T\Vert$ almost everywhere on $\set_1 \Vert T\Vert$, is \emph{Pfeffer 1 integrable} or \emph{$\PfefferOne$ integrable on $T$} if there exists a continuous additive function $F$ on $\subcspace_\leq(T)$ and for every $\epsilon>0$, there exists a gauge $\delta$ and a positive number $\tau$ such that whenever $\calP$ is a $\delta$-fine tagged family in $T$ with $\vert F(T-[\calP])\vert <\tau$, there holds:
  \begin{equation}\label{1d Pfeffer integrability}
    \vert F(T) -\sigma(f,\calP) \vert < \epsilon.
  \end{equation}
(Where $\sigma(f,\calP)$ denotes the Riemann sum $\sum_{(x,S)\in\calP} f(x)\mass(S)$.)
\end{Definition}
$F(T)$ is also the $\PfefferOne$ integral of $f$ on $T$ and we sometimes denote it $(\PfefferOne)\int_T f$.
\begin{Question}
Is it equivalent to ask that each families be surbordinate to some decomposition? This is not clear because a piece of $T$ can very well \textit{not} be a piece of any decomposition (see Figure \ref{fig:notinadecomposition}).

According to Example \ref{example two curves}, it is not sufficient to be integrable on all elements of one given decomposition to be integrable on the whole current. However, suppose $f$ is integrable on each piece for two decompositions, is the integral the same?
\end{Question}

We list the main basic properties of the integral. The proofs of the two first ones use elementary comparisons  and the fact that given two gauges $\delta_1$ and $\delta_2$, the minimum of the two is a gauge and that if $\calP$ is a $\min(\delta_1,\delta_2)$-fine family, it is also $\delta_1$ and $\delta_2$-fine. Similarly, if $\tau_1\leq \tau_2$ and $\calP$ is $(G, \tau_1)$-full in $T$, then it is $(G,\tau_2)$-full.
\begin{Proposition}
\label{1d Pfeffer linear}
The space of $\PfefferOne$ integrable functions on $T$ is a linear space and the integral: $f\mapsto I(f,T)$ is linear on this space. Furthermore, if $f\leq g$ and $f$ and $g$ are $\PfefferOne$ integrable on $T$, then $(\PfefferOne)\int_T f \leq (\PfefferOne)\int_T g$.
\end{Proposition}
\begin{Proposition}[Cauchy criterion]\label{prop 1d Pfeffer Cauchy}
  A function $f$ is $\PfefferOne$ integrable on $T$ if and only if there is a continuous nonnegative subadditive function $G$ on $\subcspace_\leq(T)$ and for every $\epsilon>0$, there exists a gauge $\delta$ and a positive $\tau$ such that for any two $\delta$-fine $(G,\tau)$ full families $\calP_1$ and $\calP_2$,
  \begin{equation}\label{eq 1d Pfeffer Cauchy}
    \vert \sigma(f,\calP_1)-\sigma(f,\calP_2) \vert < \epsilon.
  \end{equation}
\end{Proposition}


%

\begin{Proposition}\label{1d current indefinite integral}
  Let $f$ be $\PfefferOne$ integrable on the current $T\in \Iirn$. For all $S\in \subcspace_\leq(T)$, $f$ is $\PfefferOne$ integrable on $S$ and $T-S$ and $I(f,S)+ I(f,T-S) = I(f,T)$.
\end{Proposition}
\begin{proof}
  Let $G$ be a continuous nonnegative subadditive function on $\subcspace_\leq(T)$ associated to the integrability of $f$ on $T$. Fix $S\in \subcspace_\leq(T)$, notice first that $G\hel \subcspace_\leq(S)$ and $G\hel \subcspace_\leq (T-S)$ are also nonnegative continuous and subadditive. Given $\epsilon>0$. Choose a gauge $\delta$ on $\set_1\Vert T\Vert$ and a positive $\tau$  associated to $\epsilon/2$ in the definition of integrability of $f$. $\delta \hel \set_1\Vert T-S\Vert$ is a gauge on $\set_1\Vert T-S\Vert$ , so by Lemma \ref{1dcurrent HC}, there exists a $\delta$ fine $(G\hel \subcspace_\leq(T-S),\tau/2)$ full tagged family $\calP$ in $T-S$. Now given two $\delta$ fine $(G\hel \subcspace_\leq(S),\tau/2)$ full families in $S$: $\calP_1$ and $\calP_2$, we define the concatenations $\calP \cup \calP_1$ and $\calP \cup \calP_2$.
Since $[\calP] \in \subcspace_\leq(T-S)$ and $[\calP_1],[\calP_2] \in \subcspace_\leq(S)$, we have $[\calP \cup \calP_1], [\calP \cup \calP_2] \in \subcspace_\leq(T)$ so the concatenations are families in $\subcspace_\leq(T)$. They are also $\delta$ fine and for $j=1,2,$ 
\[
G(T- [\calP\cup \calP_j]) = G(T-S- [\calP]+S -\calP_j])< G(T-S -[\calP] ) + G(S- \calP_j) < \tau
\]
by subadditivity of $G$ and definition of $\calP$ an $\calP_j$. Therefore, by Proposition \ref{prop 1d Pfeffer Cauchy}
  \[
      \vert \sigma (f, \calP \cup \calP_1) - \sigma(f, \calP\cup \calP_2) \vert = \vert \sigma (f, \calP_1) - \sigma(f, \calP_2) \vert< \epsilon.
    \]
    Thus, since $\epsilon$,  $\calP_1$ and $\calP_2$ are arbitrary one can apply the Cauchy Criterion Lemma \ref{prop 1d Pfeffer Cauchy} to $S$, this proves that $f$ is $\PfefferOne$ integrable on $S$. By a similar argument $f$ is $\PfefferOne$ integrable on $T-S$. Therefore for $\epsilon>0$, choosing a gauge $\delta$ and a positive $\tau$ adapted to the integrability of $f$ on $T$, $S$ and $T-S$ at the same time, yields for $\delta$-fine $(G,T-S,\tau/2)$ and $(G,S,\tau/2)$ full families $\calP$ and $\calP'$ in $T-S$ and $S$ respectively
    \begin{multline*}
      \vert  I(f,T) - (I(f,T-S) + I(f,S))\vert \\
      \leq \vert I(f,T) - \sigma(f,\calP\cup\calP') \vert + \vert I(f,T-S) - \sigma(f,\calP) \vert  + \vert I(f,S) - \sigma(f,\calP')\vert \\
      <3\epsilon,
    \end{multline*}
   because $\calP\cup \calP'$ is a $\delta$ fine $(G,T,\tau )$ full tagged family in $T$. As $\epsilon$ is as small as we want, this concludes the proof.
 \end{proof}
 This allows us to define a function $F$ on $\subcspace_\leq(T)$ by $S\mapsto I(f,S)$, called the \emph{indefinite integral of $f$ (on $T$)}.
\begin{Proposition}\label{1d current F add and continuous}
   The indefinite integral $F$ of $f$ defined above is additive and continuous on $\subcspace_\leq(T)$.
\end{Proposition}
\begin{proof}
  For the additivity: Let $S_1$ and $S_2$ be two pieces of $T$ such that $S_1+S_2\in \subcspace_\leq(T)$. Clearly $S_1$ and $S_2$ are pieces of $S_1+S_2$, so it suffices to apply Proposition \ref{1d current indefinite integral} to see that $F(S_1) + F(S_2 ) = F(S_1+S_2)$.

  For the continuity: If $(S_j)_j$ is a sequence of pieces of $T$ converging to $0\in \subcspace_w(T)$ with $\sup_j \mass(\partial S_j) <\infty$.
  We want to show that $F(S_j)= I(f,S_j) \to 0$ as $j$ tends to infinity. By additivity, it is equivalent to show that $I(f,T-S_j) \to I(f,T)$. For $\epsilon>0$ choose a gauge $\delta$ and a positive $\tau$ associated to the integrability of $f$ on $T$. As seen above, for all $j$, $\delta$ and $\tau/2$,  are associated to $2\epsilon$ for the integrability of $f$ on $T-S_j$. Let $\calP$ be a $\delta$ fine $(G,T-S_j, \tau/2)$ full tagged family in $T-S_j$, it satisfies
  \[
    \vert \sigma(f,\calP ) - F(T-S_j) \vert <2 \epsilon.
  \]
  By continuity of $G$, if $j$ is large enough, we can suppose $G(T-[\calP]) \leq G(T-S_j) + G(S_j-[\calP]) < \tau/2 + \tau/2$, so $\calP$ is $(G,T,\tau)$ full and
  \[
     \vert \sigma(f,\calP ) - F(T) \vert < \epsilon.
   \]
   Therefore, for large enough $j$, $\vert F(T)-F(T-S_j)\vert <3\epsilon$ and we conclude that $F(S_j)\to 0$ as $j$ tends to infinity. This proves that $F$ is continuous on $\subcspace_\leq(T)$.
\end{proof}
\begin{Theorem}[Saks-Henstock Lemma]
  $f$ is $\PfefferOne$ integrable on $T$ if and only if there exists a continuous additive function $F$ on $\subcspace_\leq(T)$ satisfying: For all $\epsilon>0$, there exists a gauge $\delta$ on $\set_1\Vert T\Vert$ such that whenever $\calP$ is a $\delta$-fine tagged family in $T$:
  \begin{equation}
    \label{eq 1d current SaksHenstock}
    \sum_{(S,x)\in \calP} \vert F(S) - f(x)\mass(S)\vert < \epsilon.
  \end{equation}
\end{Theorem}
\begin{proof}
  If the second condition in the statement is satisfied, it is straightforward to prove that $f$ is $\PfefferOne$ integrable on $T$, with integral $I(f,T) = F(T)$ and the ``control function'' $G = \vert F\vert$, indeed for $\epsilon>0$, if $\delta$ is a gauge on $T$ associated to $\epsilon/2$ in the statement of the theorem and $\calP$ is a $\delta$-fine, $(G,\epsilon/2)$ full tagged family in $T$
  \[
     \vert F(T) - \sigma(f,\calP) \vert \leq \left  \vert F(T) - \sum_{(S,x)\in \calP} F(S)\right \vert + \sum_{(S,x)\in \calP} \vert F(S) - f(x)\mass(S)\vert < \epsilon.
  \]
Similarly, one proves that $F$ is the indefinite integral of $f$.

Conversely, suppose $f$ is $\PfefferOne$ integrable on $T$. The proof is very similar to the case of Henstock Kurzweil integration. Suppose that $f$ is $\PfefferOne$ integrable on $T$ and for $\epsilon>0$, fix a positive number $\tau< \epsilon/4$ and a gauge $\delta$ on $\set_1\Vert T\Vert$ such that whenever $\calP$ is a $\delta$ fine $(\vert F\vert, \tau)$ full tagged family in $T$,
\[
   \vert \sigma(f,\calP)- F(T)\vert < \dfrac{\epsilon}{4}.
 \]
 Let $\calP$ be a $\delta$ fine tagged family in $T$, without any hypothesis on $\vert F(T-[\calP])\vert$. Notice first that since $T-[\calP]$ is an integral current, there exists a $\delta$ fine, $(\vert F\vert,\tau)$ full tagged family $\calQ$ in $T-[\calP]$, which implies that $\calP\cup \calQ$ is a $\delta$ fine $(\vert F\vert, \tau/2)$ full tagged family in $T$ and
 \begin{equation}\label{eq saks henstock proof}
    \sum_{(S,x)\in \calP} \vert F(S) - f(x)\mass(S)\vert \leq  \sum_{(S,x)\in \calP\cup \calQ} \vert F(S) - f(x)\mass(S)\vert.
  \end{equation}
  Therefore, it is enough to prove that \eqref{eq 1d current SaksHenstock} holds for $(\vert F\vert, \tau)$ full families in $T$ and we suppose that $\calP$ is $(\vert F\vert, \tau)$ full. We can write $\calP = \{ (S_1,x_1),\dots, (S_p,x_p) \}$ and, reordering, assume that for some $k_0\leq p$, if $1\leq j \leq k_0$, $\vert F(S_j)-f(x_j)\mass(S_j)\vert \geq 0$ whereas for $k_0+1\leq j\leq p$, $\vert \vert F(S_j)-f(x_j)\mass(S_j)\vert < 0$. For $j=1,\dots,p$ use the $\PfefferOne$ integrability of $f$ on $S_j$ to define a $\delta$ fine, $(F\hel S_j,\tau/p)$ full tagged family $\calP_j$ such that $\vert \sigma(f,\calP_j) - F(S_j)\vert < \epsilon/(2p)$.  Consider the families
  \begin{align*}
    \calP^+ &\defeq \{(S_1,x_1), \dots, (S_{k_0}, x_{k_0})\} \cup \calP_{k_0+1} \cup \dots \cup \calP_p, \text{ and}\\
    \calP^- &\defeq  \calP_{1} \cup \dots \cup \calP_{k_0} \cup \{(S_1,x_1), \dots, (S_{k_0}, x_{k_0})\}.
  \end{align*}
  $\calP^+$ and $\calP^-$ are both $\delta$ fine, $(\vert F\vert, \tau)$ full families in $T$, therefore \eqref{eq saks henstock proof} holds for both, furthermore there holds
  \begin{multline*}
    \sum_{j=1}^{k_0} \vert F(S_j) - f(x_j)\mass(S_j)\vert = \left \vert \sum_{j=1}^{k_0} F(S_j) - f(x_j)\mass(S_j)\right \vert \\
    \leq \left \vert \sigma(f,\calP^+)-F(T)\right \vert + \sum_{j=k_0+1}^p  \vert \sigma(f,\calP_j) - F(S_j)\vert 
    \leq \dfrac{\epsilon}{4} + \dfrac{(p-k_0) \epsilon}{2p}
  \end{multline*}
    and symmetrically
    \begin{multline*}
      \sum_{j=k_0+1}^{p} \vert F(S_j) - f(x_j)\mass(S_j)\vert = \left \vert \sum_{j=k_0+1}^{p} F(S_j) - f(x_j)\mass(S_j))\right \vert \\
    \leq \left \vert \sigma(f,\calP^-)-F(T)\right \vert + \sum_{j=1}^{k_0}  \vert \sigma(f,\calP_j) - F(S_j)\vert 
    \leq \dfrac{\epsilon}{4} + \dfrac{k_0 \epsilon}{2p}.
  \end{multline*}  
  Combining the two inequalities above yields
  \[
       \sum_{j=1}^{p} \vert F(S_j) - f(x_j)\mass(S_j)\vert < \epsilon.
  \]
\end{proof}
\begin{Proposition}
  If $f$ is $\PfefferOne$ integrable on $T$, then given any decomposition  $T= T_1+T_2+\dots$, $f$ is $\PfefferOne$ integrable on $T_j$ for all $j$ with  $I(f,T)= \sum_j I(f,T_j)$. In fact $f\circ \gamma_j$ is HK integrable on $[\,0,\mass(T_j)\,]$.
\end{Proposition}
\begin{proof}
The first part of the statement is clear. For the second part, it suffices to notice that $\sum_{j=1} ^k T_j \to T$ as $k$ goes to infinity with $\sup_k \mass(\partial (\sum_{j=1}^k T_j)) \leq \mass(\partial T)$ for all $k$. By continuity of the indefinite integral $F$ of $f$ on $T$, $\sum_{j=1}^k F(T_j) \to \sum_{j=1}^\infty F(T_j)=F(T)$.
\end{proof}

\begin{Proposition}\label{Lebesgue to 1d Pfeffer}
If $f$ is defined almost everywhere in $\set_1 \Vert T\Vert$ and Lebesgue integrable with respect to $\Vert T\Vert$, then $f$ is $\PfefferOne$ integrable on $T$. As a consequence, the integral of a $\PfefferOne$ integrable function does not depend on its values on a $\Vert T\Vert$ null set.
\end{Proposition}
\begin{proof}
Let $f$ be Lebesgue integrable with respect to $\Vert T\Vert$, extend $f$ by $0$ so that it is defined everywhere in $\spt T$. Fix $\epsilon>0$. By the Vitali Caratheodory Theorem (see \cite[2.24]{RudinReal}, there exists two functions $g$ and $h$ with $g\leq f\leq h$ almost everywhere, $(\calL)\int(h-g) \dd \Vert T\Vert <\epsilon$ and $g$ and $h$ are respectively upper and lower semi-continuous. By upper (respectively lower) semi continuity of $g$ (respectively $h$), for each $x\in \spt T$, there exists $\delta(x)>0$ such that whenever $y \in \spt T\cap \Ball(x,\delta(x))$,
 \[
     g(y)\leq f(x) -\epsilon \text{ (and respectively } h(y) \geq f(x)-\epsilon \text{).}
\]
(Note that $\delta(x)$ can be chosen for $g$ and $h$ at the same time for each $x$.)  Suppose that $\calP$ is a $\delta$ fine tagged family in $T$, with $\mass(T-[\calP] ) <\epsilon$, 
\[
   (\calL) \int g \dd \Vert [\calP]\Vert - \epsilon \mass([\calP]) \leq \sigma(f,\calP) \leq (\calL) \int g \dd \Vert [\calP]\Vert + \epsilon \mass([\calP]).
\]
If $\calP_1$ and $\calP_2$ are two such families, there holds
\[
   \left \vert \sigma(f,\calP_1) -\sigma(f,\calP_2) \right \vert \leq (\calL) \int (h-g) \dd \Vert T\Vert +2 \epsilon \mass(T).
\]
As $\epsilon$ is arbitrary, we can use \cref{prop 1d Pfeffer Cauchy} to prove that $f$ is $\PfefferOne$ integrable. The $\PfefferOne$ integral of $f$ coincides with its Lebesgue integral. Indeed, choose a sequence $(\calP_j)_j$ of $\delta$ fine families in $T$ with  $\mass(T-[\calP_j])\to 0$,
\[
    (\calL)\int g\dd \Vert [\calP_j]\Vert  \to (\calL)\int g\dd \Vert T\Vert
\]
and the same holds for $h$.

In particular, if $f$ is $\PfefferOne$ integrable on $T$ and $g$ is equal to $f$, $\Vert T\Vert$ almost everywhere, then $g-f$ is equal to zero $\Vert T\Vert$ almost everywhere and is therefore Lebesgue integrable with respect to $\Vert T\Vert$, thus $\PfefferOne$ integrable on $T$ and $g = (g-f) + f$ is also $\PfefferOne$ integrable with the same integral (and indefinite integral) as $f$.
\end{proof}
\begin{Proposition}
 If $f$ is $\PfefferOne$ integrable on $T$, then its indefinite integral $F$ is $AC_*$.
\end{Proposition}
\begin{proof}
Let $N$ be a $\Vert T\Vert$ null set. By the Saks-Henstock Lemma, for $\epsilon>0$, there exists a gauge $\delta$ on $\set_1 \Vert T\Vert$ such that
\[
   \sum_{(S,x)\in \calP} \vert F(S)-f(x) \mass(S)\vert  <\epsilon,
\]
for every $\delta$ fine tagged family $\calP$ in $T$.  As $F$ does not depend on the value of $f$ on $N$, we can suppose that $f(x) = 0$ for all $x\in N$. If $\calP$ is anchored in $N$, we have
\[
   \vert F([\calP]) \vert \leq \sum_{(S,x)\in \calP} \vert F(S) \vert < \epsilon,
\]
which proves that $F$ is $AC_*$ on $T$.
\end{proof}
  \begin{Proposition}\label{1d measurable}
    If $f$ is $\PfefferOne$ integrable, then it is $\Vert T \Vert$ measurable.
  \end{Proposition}
  \begin{proof}
    Consider a decomposition of $T$: $T_1,T_2,\dots$ and a representative of $f$. $f$ is $\PfefferOne$ integrable on each $T_k\eqdef \lseg \gamma_{k}\rseg$ and therefore, $f\circ \gamma_k$ is HK integrable on $[\,0,\mass(T_k)\,]$ and thus Lebesgue measurable. Thus $f$ is $\Vert T_k\Vert$ measurable, and also, $f_k \defeq f\hel \spt \Vert T_k\Vert$ is $\Vert T\Vert$ measurable. Consider the function $\tilde f:\, x\mapsto \sup_k f_k(x)$. $\tilde f $ is $\Vert T\Vert$ measurable as a pointwise supremum of measurable functions. The function $f-\tilde f$ is equal to zero at each point of $\spt T_1 \cup \spt T_2 \cup \dots \subseteq \set_1 \Vert T\Vert$. By definition of decomposition of currents, $\Vert T\Vert = \sum_{k=1}^\infty \Vert T_k\Vert$ and as indecomposable currents correspond to simple Lipschitz curve with integral multiplicity, for all $k$, $\spt  T_k = \set_1 \Vert T_k\Vert$, therefore 
\[
    \Vert T\Vert(\Rn\backslash \bigcup_{k=1}^\infty \spt  T_k) = 0,
\]
thus $f = \tilde f$, $\Vert T\Vert$ almost everywhere. This proves that $f$ is $\Vert T\Vert$ measurable.
\end{proof}

\begin{Proposition}\label{1d Pfeffer nonnegative Lebesgue}
  Conversely to Proposition \ref{Lebesgue to 1d Pfeffer}, a function $f$ is Lebesgue integrable with respect to $\Vert T\Vert$ if and only if $f$ and $\vert f\vert$ are $\PfefferOne$ integrable on $T$.
\end{Proposition}
\begin{proof}
Without loss of generality, we can suppose that $f$ is nonnegative and $\PfefferOne$ integrable, we also fix a representative of $f$ with respect to $\Vert T\Vert$. It suffices to show that $f$ is Lebesgue integrable with respect to $\Vert T\Vert$. For $k=1,2,\dots$, consider the function $f_k \defeq f \ind_{\{x, f(x) \leq k\}}$. Since $f$ is $\Vert T\Vert$ measurable by Proposition \ref{1d measurable}, $f_k$ is $\Vert T\Vert$ measurable and bounded and thus Lebesgue integrable with respect to $\Vert T\Vert$ (which is a finite measure).
the sequence $f_k$ is nondecreasing and converges pointwise to $f$. Furthermore the sequence $\left ((\calL) \int f_k \dd \Vert T\Vert\right)_k = \left ( (\PfefferOne) \int_T f_k \right)_k$ is bounded from above by $(\PfefferOne)\int_T f$. By the Lebesgue Monotone Convergence Theorem, $f$ is Lebesgue integrable with respect to $\Vert T\Vert$.
\end{proof}
  \begin{Theorem}[Monotone Convergence Theorem for the $\PfefferOne$ integral.]\label{1d monotone convergence}
    Suppose that  $(f_k)_{k=1,2,\dots}$ is a $\Vert T\Vert$ almost everywhere nondecreasing sequence of $\PfefferOne$ integrable functions on $T$. If there exists $f: \set_1 \Vert T\Vert\to \R$ such that $f_k(x)$ converges to $f(x)$ $\Vert T\Vert$ almost everywhere and if furthermore, the sequence of integral: $(\PfefferOne) \int_T f_k$ for $k=1,2,\dots$ is bounded from above. Then $f$ is $\PfefferOne$ integrable on $T$ with 
\[ 
 (\PfefferOne)\int_T f = \lim (\PfefferOne) \int_T f_k.
\]
  \end{Theorem}
We give a proof which does not rely on the measurability of $f$ or on Lebesgue integration results, but relies only on gauge integration techniques.
\begin{proof}
Since the $\PfefferOne$ integral of a function does not depend on its values in a $\Vert T \Vert$ null set, we can suppose that $f_k$ converges pointwise to $f$ everywhere and that for all $x\in \set_1\Vert T \Vert$, the sequence $(f_k(x))_k$ is nondecreasing. Up to substracting $f_1$, we can also suppose that all the $f_k$ are nonnegative (by linearity of the integral). For $k=1,2,\dots$, let $F_k$ be the indefinite $\PfefferOne$ integral of $f_k$ on $T$ it is nonnegative. Notice also that for all $S\in \subcspace_\leq(T)$, and for $k \leq k'$, $F_k(S) \leq F_{k'}(S)$ by the last part of \cref{1d Pfeffer linear}.
Since $F_k(T)$ is bounded from above, it converges to a limit $F(T)$, similarly we can define $F(S)$ for any $S\in \subcspace_\leq(T)$ as both $(F_k(T-S))_k$ and $(F_k(S))$ are nondecreasing sequences bounded from above by $F(T)\geq F_k(S) + F_k(T-S)$. $F$ is nonnegative. The function $F$ on $\subcspace_\leq(T)$ is also additive, indeed, suppose $S$, $S'$ and $S+S'$ are in $\subcspace_\leq(T)$, we have
\[
   F(S+S') = \lim_{k\to \infty} F_k(S+S')  = \lim_{k\to \infty} (F_k(S) + F_k(S')) = F(S)+F(S').
\]
Let us now prove that $F$ is continuous. Fix sequence $(S_j)_j$ in $\subcspace_\leq(T)$ with $\sup_j \mass(\partial S_j) < \infty$ and $\flatn(S_j) \to 0$.  For each $k$, the sequence $(F_k(S_j))_j$ goes to $0$ as $j$ goes to $\infty$ and similarly $F_k(T-S_j)\to F_k(T)$ as $j\to \infty$. Thus, since for all $k$ and $j$, $F(T) \geq F(T-S_j) \geq F_k(T-S_j)$, given $\epsilon>0$ there exists $k_0$ such that for all $k\geq k_0$, $F_k(T) \geq F(T)-\epsilon/2$.

There exists also $j_0$ such that for all $j\geq j_0$, $F_{k_0}(T-S_j) \geq F_{k_0}(T)-\epsilon/2$. This implies that for all $j\geq j_0$ and all $k\geq k_0$, 
\[
   F(T) \geq F(T-S_j) \geq F_k(T-S_j) \geq F_{k_0}(T-S_j)\geq F_{k_0}(T)-\dfrac{\epsilon}{2} \geq F(T)-\epsilon.
\]

Thus $F$ is nonnegative, additive and continuous on $\subcspace_\leq(T)$.
Since $F(S) \geq F_k(S)$ for all $k$, if $\calP$ is an $(F,\tau)$ full tagged family in $T$ for some $\tau>0$, $\calP$ is also $(F_k,\tau)$ full for all $k$.

From now on the argument follows the method of \cite[4.42]{Moonens2017integration}.
Fix $\epsilon>0$, there exists $l$ such that for all $k\geq l$, $F(T)- F_k(T) <\epsilon/4$.
For each $k\geq l$, fix a gauge $\delta'_k$ on $\set_1 \Vert T\Vert$ such that for all $\delta'_k$ fine, ($\vert F_k\vert$, $\epsilon/4)$ full families $\calP$ int $T$,
\[
    \sum_{(x,S)\in \calP} \vert F_k(S)-f_k(x)\mass(S) \vert < \dfrac{\epsilon}{4^{k+2}}.
\]
Define a new series of gauges $(\delta_k)_k$ such that for $x\in \set_1\Vert T\Vert$,
\[
    \delta_k(x) \defeq \min_{1\leq j\leq k} \delta'_k(x).
\]
Note that $\delta_k$ is indeed a gauge, as a finite union of countable sets is countable.
For each $x\in \set_1 \Vert T\Vert$, fix $l(x)\geq l$ so that $0\leq f(x)- f_k(x)<\epsilon /(4\mass(T))$ whenever $k\geq l(x)$. And let $\delta(x)\defeq \delta_{l(x)}(x)$ be a gauge on $\set_1\Vert T\Vert$. To check that the zero set of $\delta$ is countable, notice that it is contained in the countable union of the zero sets of the gauges $\delta_k'$.

Let $\calP$ be a $\delta$ fine, $(F,\epsilon/4)$ full tagged family in $T$. It is also $(F_k,\epsilon/4)$ full, as we said above.
Let $l'$ be the maximum of the indices $l(x)$ over $(x,S)\in \calP$. For $l\leq k\leq l'$ let $\calP_k$ be the subfamily of $\calP$ consisting of all the $(x,S)\in \calP$ with $l(x) = k$. We can write
\begin{multline*}
     \sigma(f,\calP) -F(T) = \sum_{k=l}^{l'} \sigma(f,\calP_k) -F(T)\\
= \sum_{k=l}^{l'} (\sigma(f,\calP_k)-\sigma(f_j,\calP_k)) + \sum_{k=l}^{l'} \left (\sigma(f_k,\calP_k)-F_k([\calP_k]) \right) + \sum_{k=l}^{l'} F_k([\calP_k]) - F(T).
\end{multline*}
To control the first term, by the choice of $l(x)$, for all $k$ we have 
\[
    0\leq \sigma(f,\calP_k) -\sigma(f_k \calP_k) < \dfrac{\mass([\calP_k])}{\mass(T)}\dfrac{\epsilon}{4}.
\]
Sum over $k= l,\dots,l'$ to obtain
\[
    0\leq \sum_{k=l}^{l'} \sigma(f,\calP_k) -\sigma(f_k \calP_k) < \dfrac{\mass([\calP])}{\mass(T)}\dfrac{\epsilon}{4} \leq \dfrac{\epsilon}{4}.
\]
For the second term, for any $k$ by the Saks-Henstock Lemma applied to $f_k$ and $\calP_k$ we have
\[
    \left \vert \sigma(f_k,\calP_k)-F_k([\calP_k])  \right\vert \leq \dfrac{\epsilon}{4^{k+2}}.
\]
Which can be summed to get
\[
    \sum_{k=l}^{l'} \left \vert \sigma(f_k,\calP_k)-F_k([\calP_k])  \right\vert \leq \dfrac{\epsilon}{4}.
\]
Finally, for the third term, notice that for all $k\geq l$
\[
   F_k([\calP_k])  \geq F_l([\calP_k]).
\]
Summing over $k$ yields
\[
    F(T)\geq F([\calP]) \geq \sum_{k=l}^{l'} F_k([\calP_k]) \geq F_l ([\calP])\geq F_l(T)-\dfrac{\epsilon}{4}\geq F(T) - \epsilon/2,
\]
as $\calP$ is $(F_l,\epsilon/4)$ full in $T$.
Combining the three above estimates we get
\[
   \vert \sigma(f,\calP)-F(T)\vert < \epsilon,
\]
which proves that $f$ has $\PfefferOne$ integral $F(T)$ on $T$. By the same reasonning one can prove that $f$ is $\PfefferOne$ integrable on $S\in \subcspace_\leq(T)$ with integral $F(S)$, thus $F$ is the indefinite integral of $f$ on $T$.
\end{proof}
\subsection{Fundamental Theorem of Calculus for the $\PfefferOne$ integral} 

\begin{Proposition}\label{1d FTI  ACstar charge}
  If $F$ is a continuous additive function on $\subcspace_\leq(T)$ which is $AC_*$ and derivable $\Vert T\Vert$ almost everywhere, then $x\mapsto \dcharge_T F(x)$ is $\PfefferOne$ integrable on $T$ with indefinite integral $F$.
\end{Proposition}
\begin{proof}
Let $N$ be the set of non derivability points of $F$ in $\set_1\Vert T\Vert$. Let $f$ be the function defined on $\set_1\Vert T\Vert$ by $f(x) = 0$ if $x\in N$ and $f(x) = \dcharge_T F(x)$ otherwise.  For $\epsilon>0$, let $\delta$ be a gauge on $\set_1\Vert T\Vert$ such that whenever  $\calP$ is a $\delta$-fine tagged family in $T$ anchored in $N$, $\vert F([\calP]) \vert <\epsilon$ and for all $x\in \set_1 \Vert T\Vert \backslash N$, $\delta(x)$ is a positive number such that for all $S\in \subcspace_\leq(T,x,\delta(x))$
\[
    \left \vert F(S) -f(x) \mass(S) \right \vert < \epsilon \mass(S).
\]
If $\calP$ is a $\delta$-fine tagged family in $T$ with $\vert F(T-[\calP])\vert <\epsilon$, let $\calP_N$ be the subfamily of $\calP$ containing all the pairs $(S,x)\in \calP$ with $x\in N$. There holds
\begin{align*}
    \vert F(T)-\sigma(f,\calP) \vert\hspace{-2cm} \\
&\leq \vert F(T-[\calP])\vert +\vert F([\calP_N])\vert + \sum_{(S,x)\in \calP, x\notin N} \left \vert  F(S) -f(x) \mass(S) \right \vert \\
&< 3\epsilon
  \end{align*}
  Thus $f$ is $\PfefferOne$ integrable in $T$ with  $I(f,T) = F(T)$. Since $F\vert_{\subcspace_\leq(S)}$ satisfies the hypothesis of the theorem for any $S\in \subcspace_\leq(T)$, $I(f,S)=F(S)$ and $F$ is the indefinite integral of $\dcharge_T F$ on $T$. 
\end{proof}
\begin{Proposition}\label{1d prop almost diff and ACstar}
If $u$ is a continuous function on $\spt T$ which is differentiable $\Vert T\Vert$ almost everywhere and $\Theta_u$ is $AC_*$, then the function
\[
   x\mapsto \dcharge_T \Theta_u (x) = \langle \dD u(x),\vect T(x)\rangle
\]
 is $\PfefferOne$ integrable on $T$ with indefinite integral $\Theta_u$.
\end{Proposition}
\begin{proof}
Using Proposition \ref{1d FTI ACstar charge} it suffices to prove that the set
\[
  \{x, \Theta_u \text{ is not derivable at } x\}\cup \{x, \dcharge_T\Theta_u(x) \neq \langle \dD u(x),\vect T(x) \rangle\}
\]
is $\Vert T\Vert$ negligible.
As $u$ is differentiable $\Vert T\Vert$ almost everywhere, by Theorem \ref{1d thm derivability} \ref{1d u derivable} this reduces to proving that the set of points $x$ at which $\vect{T}$  has a $\Vert T\Vert$ approximately continuous representative is $\Vert T\Vert$ negligible.
\begin{Claim}\label{vect T approx cont}  
The function $x \mapsto \vect T$ is $\Vert T\Vert$ approximately continuous $\Vert T\Vert$ almost everywhere, i.e. for $\Vert T\Vert$ almost every $x$, for every $\epsilon>0$, there exists $\delta>0$ such that
\[
   \Theta^{m*}(\Vert T\Vert\hel \{y, \vert \vect T(x)-\vect T(y) \vert \geq \delta \},x) <\epsilon.
\]
\end{Claim}
\begin{innerproof}
The measure $\Vert T\Vert$ in $\Rn$ is finite and Borel regular, therefore the Besicovitch Covering Theorem (see \cite[Theorem 2.7]{Mattila1995book}) holds for $\Vert T\Vert$. In the words of H. Federer \cite[2.8.9, 2.8.18]{FedererGMT}, the ambient space $\Rn$ is directionally limited and the collection of balls
\[
   \{(x,\Ball(x,r)\, \vert \, x\in \Rn, r>0\},
\]
 forms a Vitali relation for the measure $\Vert T\Vert$. Furthermore, the function $\vect T: \set_1\Vert T\Vert \to \Lambda_1(\Rn)$ is $\Vert T\Vert$ measurable. Thus, by\cite[2.9.13]{FedererGMT}, the vector function $\vect T$ is $\Vert T\Vert$ approximately continuous $\Vert T\Vert$ almost everywhere. 
\end{innerproof}
\end{proof}
We can finally restate and prove our main result:
\hkpfti*
\begin{proof}
Let $\Theta_u$ be the function on $\subcspace_\leq(T)$ associated to the variations of $u$. By Proposition \ref{1d thm derivability}\ref{1d u pointwise lip}, $\Theta_u$ is almost derivable at all points of $\set_1\Vert T\Vert$ except for a countable set. By Theorem \ref{1d prop almost diff ACstar}, $\Theta_u$ is $AC_*$. By \cref{1d thm derivability} \ref{1d u derivable},  $\Theta_u$ is derivable $\Vert T\Vert$ almost everywhere along $T$ with derivative equal to $\langle \dD u (x), \vect T(x)\rangle$. Use propositions \ref{1d FTI ACstar charge} and \ref{1d prop almost diff and ACstar} to conclude.
\end{proof}


\bibliographystyle{plain}
\bibliography{../Bibliographie/RefsmathsAJ}
\end{document}